\documentclass{birkau}
\usepackage{tikz}
\usepackage{amsmath,amssymb,url}
\usepackage{enumitem} 
\usepackage{bussproofs}
\usepackage{proof}

\usepackage[utf8]{inputenc}   
\usepackage[T1]{fontenc}      
\usepackage{lmodern}          
\usepackage{hyperref}
\newcommand{\boxg}{{\Box_{\mathfrak o}}}

\newcommand{\boxm}{{\Box_{\mathfrak p}}}

\newcommand{\diamondg}{{\Diamond_{\mathfrak o}}}

\newcommand{\diamondm}{{\Diamond_{\mathfrak p}}}

\numberwithin{equation}{section}

\theoremstyle{plain}
\newtheorem{theorem}{Theorem}[section]
\newtheorem{lemma}[theorem]{Lemma}
\newtheorem{proposition}[theorem]{Proposition}
\newtheorem{corollary}[theorem]{Corollary}
\newtheorem{example}[theorem]{Example}
\theoremstyle{definition}
\newtheorem{definition}[theorem]{Definition}

\newtheorem{remark}[theorem]{Remark}
\newtheorem*{note}{Note} 

\def\KK{\mathbb{K}}
\def\BK{\mathfrak{B}(\KK)} 
\def\PK{\mathfrak{P}(\KK)} 
\def\HK{\mathfrak{H}(\KK)} 

\begin{document}


\title[DBas: axiomatization and topological representation]{Towards a Simplified Theory of Double Boolean Algebras: Axioms and Topological Representation}




\author[P. Howlader]{Prosenjit Howlader}
\address{ Institute of Information Science\\
Academia Sinica\\Taipei,115\\Taiwan}
\email{prosen@mail.iis.sinica.edu.tw}

\corrauthor[L. Kwuida]{Leonard Kwuida}
\address{Bern University of Applied Sciences\\ School of Business\\Br\"uckenstr.~73\\ 3005 Bern\\ Switzerland}
\email{leonard.kwuida@bfh.ch}
\author[M. Behrisch]{Mike Behrisch}
\address{Technische Universit\"at Wien\\ Institut f\"ur Diskrete Mathematik und Geometrie\\Wiedner Hauptstr. 8--10, \\1040 Wien\\ Austria}
\email{behrisch@logic.at}

\author[C. J. Liau]{Churn-Jung Liau}
\address{ Institute of Information Science\\
Academia Sinica\\Taipei,115\\Taiwan}
\email{ liaucj@iis.sinica.edu.tw}





\begin{abstract}
Double Boolean algebras (dBas), introduced by Wille, are based on twenty-three identities. We present a simplified axiom system, the \emph{D-core algebra}, and prove it is equivalent to Wille’s original definition. This reduction allows improved structural results, including a refined Boolean representation theorem showing fewer conditions suffice to represent a dBa as a pair of Boolean algebras linked by adjoint maps. We generalize the glued-sum construction to possibly overlapping Boolean algebras, characterize them via a generalized order, and establish a Stone-type topological representation: every dBa is quasi-isomorphic to a dBa of clopen subsets of a Stone space. Simplified logical systems for contextual and pure dBas are developed with soundness and completeness.

\end{abstract}

\maketitle


\section{Introduction}\label{intro}

In lattice theory, a  \emph{polarity} \cite{birkhoff1940lattice} is a triplet $\KK:=(G,M, I)$ where $G$ and $M$ are sets and $I\subseteq G\times M$.  In Formal Concept Analysis (FCA), a polarity is called a \emph{formal context}~\cite{wille1982restructuring}. From a formal context, clusters can be formed in terms of formal concepts or proto-concepts, and rules can be extracted in terms of implications or association rules \cite{Lakhal2005}. These patterns are a good indication why FCA has many applications in knowledge discovery and management. To formalize these notions, the \emph{derivation operators} are defined on subsets $A\subseteq G$ and $B\subseteq M$  by\footnote{For $x\in G\cup M$, we write $x^\prime$ for $\{x\}^\prime$.}:
\begin{align*}
A^{\prime}:=\{m\in M\mid gIm, \forall g\in A \}\quad \text{ and }\quad B^{\prime}:=\{g\in G\mid gIm, \forall m\in B \}.    
\end{align*}
The pair of maps $(^\prime,^\prime)$ forms a Galois connection between the power set of $G$ and that of $M$. The composition $^{\prime\prime}$ is a closure operator, and the corresponding closed sets form a complete lattice. 

A \emph{formal concept} is a pair $(A,B)$ with $A^{\prime}=B$ and $B^{\prime}=A$, i.e., $A$ is the set of objects that satisfy all attributes of $B$, and $B$ is the set of attributes common to all objects of $A$. 
A \emph{protoconcept} is a pair $(A, B)$ with $A^{\prime\prime}=B^{\prime}$. If $(A,B)$ is a concept and $C\subseteq G$, $D\subseteq M$ such that $C^{\prime\prime}=A$ and $D^{\prime\prime}=B$, we say that $C,D$ each generates the concept $(A,B)$. Thus, $(A,B)$ is a protoconcept if $A$ and $B$ generate the same concept. Special protoconcepts of the form $(A,A^\prime)$ or $(B^\prime,B)$ are called \emph{semiconcepts}. 

We denote by $\BK$, $\HK$ and $\PK$ the set of formal concepts, semiconcepts and protoconcepts of $\KK$, respectively. It is a straightforward observation that $\BK\subseteq\HK\subseteq\PK$.  
The concept hierarchy, captured by the order relation $\le$ on $\BK$ given by
\begin{align}
    (A,B)\le (C,D):\iff A\subseteq C (\text{ or equivalently } B\supseteq D),
\end{align}
turns $\BK$ into a complete lattice~\cite{wille1982restructuring}, called \emph{concept lattice} of the formal context $\KK$. This subconcept-superconcept order generalizes to protoconcepts: 
\begin{align}
    (A,B)\sqsubseteq (C,D):\iff A\subseteq C \text{ and  } B\supseteq D.
\end{align}
 The meet and join operations of the complete lattice $(\BK,\le)$ are given by 
 \begin{align}
    (A,B)\wedge (C,D) = (A\cap C,  (B\cup D)^{\prime\prime}) = (A\cap C,  (A\cap C)^{\prime})\\
    (A,B)\vee (C,D) = ((A\cup C)^{\prime\prime},  B\cap D) = ((B\cap D)^{\prime},  B\cap D)
\end{align}
and can be extended to the set of protoconcepts. Additionally, two negation operators $\neg$ and $\lrcorner$ are defined on the set $\PK$: 
\begin{align*}
	\text{meet:} && (A,B)\sqcap(C,D) &:=(A\cap C,(A\cap C)^{\prime}),\\
    \text{join:} && (A,B)\sqcup(C,D) &:=((B\cap D)^{\prime},B\cap D),\\    
	\text{negation:} && \neg(A,B)&:=(G\setminus A,(G\setminus A)^{\prime}),\\
    \text{opposition:} && \lrcorner (A,B)&:=((M\setminus B)^{\prime},M\setminus B),\\
    \text{all:} &&  \top &:=(G,\emptyset) \text{ and } \\
    \text{nothing:} && \bot&:=(\emptyset,M).
\end{align*}
Equipped with the above defined operations 
the set $\PK$ forms an algebraic structure $\underline{\mathfrak{P}}(\mathbb{K}):=(\PK; \sqcap, \sqcup, \neg, \lrcorner, \top, \bot)$, called the \emph{algebra of protoconcepts} of $\KK$ \cite{wille}. 
The set of all semiconcepts $\HK$ forms a subalgebra of the algebra of protoconcepts,  
called the \emph{algebra of semiconcepts} of $\KK$ \cite{wille}, and 
is denoted by $\underline{\mathfrak{H}}(\mathbb{K}):=(\HK; \sqcap, \sqcup, \neg, \lrcorner, \top, \bot)$.

To abstract the algebra of protoconcepts (and of semiconcepts), double Boolean algebras 
are introduced. A double Boolean algebra (dBa) \cite{wille} is an algebra 
$(D;\sqcap, \sqcup, \neg, \lrcorner, \top, \bot)$ of type $(2,2,1,1,0,0)$ defined below. 
\begin{definition} \label{def:DBA}
	{\rm \cite{wille} An  algebra $  \mathbf{D}:= (D;  \sqcap, \sqcup,\neg,\lrcorner,\top,\bot)$ satisfying the following properties is called a  \emph{double Boolean algebra} (dBa). For any $x,y,z \in D$,\newline
		$\begin{array}{llll}
			(1a)&  (x \sqcap x ) \sqcap  y = x \sqcap  y  &
			(1b)&  (x \sqcup x)\sqcup  y = x \sqcup y \\
			(2a)& x\sqcap y = y\sqcap  x  &
			(2b)&  x \sqcup   y = y\sqcup   x  \\
			(3a)& \neg (x \sqcap  x) = \neg  x  &
			(3b)&  \lrcorner(x \sqcup   x )= \lrcorner x \\
			(4a)&  x  \sqcap (x \sqcup y)=x \sqcap  x  &
			(4b)&  x \sqcup  (x \sqcap y) = x \sqcup   x \\
			(5a)& x \sqcap  (y \vee z ) = (x\sqcap  y)\vee (x \sqcap  z) &
			(5b)&  x \sqcup  (y \wedge z) = (x \sqcup  y) \wedge  (x \sqcup  z) \\
			(6a)&  x \sqcap (x\vee y)= x \sqcap  x  &
			(6b)&  x\sqcup  (x \wedge  y) =x \sqcup   x \\
			(7a)&  \neg \neg (x \sqcap  y)= x \sqcap  y &
			(7b)&  \lrcorner\lrcorner(x \sqcup  y) = x\sqcup  y \\
			(8a)&  x  \sqcap \neg  x= \bot &
			(8b)&  x \sqcup \lrcorner x = \top  \\
			(9a)&  \neg \top = \bot  &
			(9b)& \lrcorner\bot =\top \\
			(10a)&  x \sqcap ( y \sqcap  z) = (x \sqcap  y) \sqcap  z  &
			(10b)&  x \sqcup (y \sqcup  z) = (x \sqcup  y)\sqcup  z \\
			(11a)& \neg \bot = \top \sqcap   \top  &
			(11b)& \lrcorner\top =\bot \sqcup  \bot \\
		\end{array}$
		$\begin{array}{llll}
    (12)&  (x \sqcap  x) \sqcup (x \sqcap x) = (x \sqcup x) \sqcap (x \sqcup x) &
         \end{array}$  
		where $ x\vee y := \neg(\neg x \sqcap\neg y)$ and 
		$ x \wedge y :=\lrcorner(\lrcorner x \sqcup \lrcorner y)$. 
	}
	
\end{definition}

   The operations $\sqcap$ and $\sqcup$ are called meet and join. $\neg$ and $\lrcorner$ are two negations. $\bot$ and $\top$ are called bottom and top elements. Every Boolean algebra is a  dBa by duplicating the complement, and a dBa becomes a Boolean algebra when $\neg$ and $\lrcorner$ coincide and $\neg\neg x =x$.
   For a dBa $\textbf{D}$, we use the following notations 
 $D_{\sqcap}:=\{x\in D~:~x\sqcap x=x\}$,  $D_{\sqcup}:=\{x\in D~:~x\sqcup x=x\}$ and $D_p=D_\sqcap\cup D_\sqcup$.
For $x\in D$, $x_{\sqcap}:=x\sqcap x$ and $x_{\sqcup}:=x\sqcup x$.
     
A dBa is \emph{pure} if it satisfies $x\sqcap x=x$ or $x\sqcup x=x$, i.e. $D=D_p$. A typical example is the algebra of semiconcepts, $\underline{\mathfrak{H}}(\KK)$. Note that $(\PK;\sqsubseteq)$ is a poset satisfying the condition \quad $x\sqsubseteq y \iff x\sqcap y=x\sqcap x$ and  $x\sqcup y=y\sqcup y$.
 A relation $\sqsubseteq$ is then defined on any dBa $\mathbf{D}$ by 
 \begin{align}\label{Eq:dBa quasiorder}
    x\sqsubseteq y :\iff x\sqcap y=x\sqcap x \text{ and } x\sqcup y=y\sqcup y. 
\end{align}
 and turns $(D;\sqsubseteq)$ into a quasi-ordered set. We say that $\textbf{D}$ is \emph{contextual} or \emph{regular} if the relation $\sqsubseteq$ 
is an order relation.  A contextual dBa is said to be \emph{fully contextual} if for each $y\in D_{\sqcap}$ and $x\in D_{\sqcup}$ with $y_{\sqcup}=x_{\sqcap}$, there is a unique $z\in D$ with $z_{\sqcap}=x$ and $z_{\sqcup}=y$. A typical example is the algebra of protoconcepts, $\underline{\mathfrak{P}}(\KK)$. A dBa 
is said to be \emph{trivial} if it satisfies $\top \sqcap \top =\bot\sqcup\bot$~\cite{kwuida2007prime}. Observe that if $\underline{\mathfrak{P}}(\KK)$ is a trivial dBa, then $|\BK|=1$.  

%
As the name suggests, each dBa induces two Boolean algebras: 
$\textbf{D}_{\sqcap}:=(D_{\sqcap}; \sqcap, \vee, \neg,  \bot, \top\sqcap \top)$ and
$\textbf{D}_{\sqcup}:=(D_{\sqcup};\wedge, \sqcup, \lrcorner,\bot\sqcup\bot, \top)$. 
A characterization of double Boolean algebras was presented in \cite{bmlpl}, showing that every dBa 
can be constructed from the two Boolean algebras $\textbf{D}_{\sqcap}$ and $\textbf{D}_{\sqcup}$. In \cite{kembang2023simple}, a trivial and pure dBa is characterized as linear sum of two Boolean algebras.

  Apart from the Boolean representation, several authors have studied representation theorems for dBas based on contexts. Wille~\cite{wille} constructed a \emph{standard context} for each dBa $\mathbf{D}$, denoted $\mathbb{K}(\mathbf{D}):=(\mathcal{F}_{pr}(\mathbf{D}), \mathcal{I}_{pr}(\mathbf{D}), \Delta)$, consisting of the sets $\mathcal{F}_{pr}(\mathbf{D})$ of all primary filters and  the set $\mathcal{I}_{pr}(\mathbf{D})$ of primary ideals of $\mathbf{D}$, with a relation $\Delta$ such that $F \Delta I$ iff $F \cap I \neq \emptyset$. It was shown that for every dBa $\mathbf{D}$  there is a homomorphism from $\textbf{D}$  to the algebra of \emph{protoconcepts} of $\mathbb{K}(\mathbf{D})$ \cite{wille} that preserves and reflects the quasi order (we call it a \emph{quasi-embedding}). The map becomes an embedding if $\mathbf{D}$ is \emph{pure}~\cite{BALBIANI2012260}.

The topological representation of dBas was studied in~\cite{MR4566932}, where the prime ideal theorem~\cite{howlader3,kwuida2007prime} plays a central role. To represent fully contextual and pure dBas, \emph{object-oriented protoconcepts} and \emph{semiconcepts} were introduced~\cite{howlader2018algebras,howlader2020}. A context is extended to a \emph{context on topological spaces} (CTS) by equipping the sets of objects $G$ and properties $M$ with topologies. Within this topological framework, \emph{clopen} object-oriented \emph{semiconcepts} and \emph{protoconcepts} are defined and shown to form a fully contextual dBa and a pure dBa, respectively.
   In~\cite{MR4566932}, the complement of the standard context,  $\mathbb{K}^c(\mathbf{D})$, is considered, where the relation $\Delta$ in $\mathbb{K}(\mathbf{D})$ is replaced with its complement $\nabla$, defined by $F \nabla I$ if and only if $F \cap I = \emptyset$. Then, $\mathbb{K}^c(\mathbf{D})$ is extended with topologies on primary filters and ideals to yield the CTS $\mathbb{K}_{pr}^{T}(\mathbf{D})$. It is shown that $\mathbb{K}_{pr}^{T}(\mathbf{D})$ forms a CTS 
   where both the incidence relation and its converse are continuous~\cite{guide2006infinite}, using the prime ideal theorem~\cite{kwuida2007prime,howlader3}. We call it a \emph{context on topological spaces with continuous relations} (CTSCR, see Definition~\ref{def:CTSCR topcon}). In \cite{MR4566932}, It is proved that every fully contextual  dBa is isomorphic to  the algebra of clopen object-oriented  protoconcepts of $\mathbb{K}_{pr}^{T}(\mathbf{D})$, and  if $\mathbf{D}$ is pure, it is isomorphic to the algebra of clopen object-oriented semiconcepts. This result also provides a representation for Boolean algebras as special cases of pure dBas. The topological representation of contextual and arbitrary dBas remained an open question, which we address in this article.

Beyond their algebraic analysis, dBas have also been systematically studied from the perspective of logic, with both approaches enriching our understanding of the structure. In \cite{howlader2021dbalogic}, a sequent calculus, \textbf{CDBL}, was developed for contextual dBas and subsequently extended to its modal version. In \cite{HOWLADER2023115}, a hypersequent calculus, \textbf{PDBL}, was introduced for pure dBas and similarly extended to a modal system. Relational semantics for both calculi have been studied on the basis of formal contexts. In this article, we propose an equivalent proof system for \textbf{CDBL} and another equivalent proof system for \textbf{PDBL}. The new proof systems presented here are significantly simplified. 
   
Double Boolean algebras are defined by a set of 23 equations, as presented in Definition~\ref{def:DBA}. This relatively long list of axioms may discourage potential users from exploring or applying these structures. In previous work~\cite{howladerdiscussion}, it was shown that axioms $(1a)$, $(1b)$, $(11a)$, and $(11b)$ are derivable from the other axioms of a dBa. 
   In this article, we develop a more streamlined, redundancy-free axiomatization. In particular,  we introduce the notion of a \emph{D-core} algebra. 
That is a universal algebra of type $(2,2,1,1,0,0)$ 
that satisfies conditions (2a)--(5a), (7a)--(8a), (2b)--(5b), (7b)--(8b), 
and 12, as listed in Definition~\ref{def:DBA}. Based on this definition, we establish the following result.

   \begin{itemize}
       \item A universal algebra $(D,  \sqcap, \sqcup, \neg, \lrcorner, \top, \bot)$ is a dBa if and only if it  is a D-core algebra.
   \end{itemize}
   Moreover, one can verify that D-core algebras represent the minimal set of axioms for dBas. As a demonstration, we construct a non-dBa that satisfies all the axioms of a D-core algebra except (7a) and (7b).

   The new, reduced definition of dBas helps us better understand structural results. We prove a stronger version of Theorem~\ref{dbaboolean} by weakening its original conditions. The original result, proved in \cite{bmlpl}, describes how a dBa can be constructed from two Boolean algebras. This new Boolean representation theorem is used to prove further structural results for the D-core algebra. 
We define the \emph{generalized glued sum} of two posets, denoted by $P \oplus_{g} Q$, in the same way as the ordinary glued sum of posets, except that we do not require $P$ and $Q$ to be disjoint. 
Unlike the usual glued sum, the generalized glued sum lacks the antisymmetry property; hence, $P \oplus_{g} Q$ is, in general, only a quasi-order. 
Next, this construction is extended to Boolean algebras. 
The generalized glued sum combines two Boolean algebras by uniting their elements, preserving each algebra’s internal operations, and defining a global order based on how elements interact under these combined operations.  

Next, we illustrate how to define the generalized glued sum of two Boolean algebras, {\bf P} and {\bf Q} with domains $P$ and $Q$, respectively. We consider the following maps: 
$r$ maps elements of $P \cup Q$ to $P$, sending elements of $P$ to themselves and elements outside $P$ to the top element of $P$; 
$e$ embeds elements of $P$ into $P \cup Q$. Similarly, $r'$ maps elements of $P \cup Q$ to $Q$, sending elements of $Q$ to themselves and elements outside $Q$ to the bottom element of $Q$; $e'$ embeds elements of $Q$ into $P \cup Q$. Using these maps\footnote{Intuitively, the pairs of maps $(r,e)$ and $(r',e')$ encode the maps $D\to D_\sqcap, x\mapsto x\sqcap x$, $D_\sqcap\to D, x_\sqcap \mapsto x_\sqcap$ and $D\to D_\sqcup, x\mapsto x\sqcup x$, $D_\sqcup\to D, x_\sqcup \mapsto x_\sqcup$, respectively, whenever $D$ is a dBa. i.e $P=D_\sqcap$ and $Q=D_\sqcup$.}, we define a universal algebra
\[
\mathbf{P+Q} := (P \cup Q, \sqcap, \sqcup, \neg, \lrcorner, e'(\top_q), e(\bot_p)),
\]
where the operations are given by 
\begin{align}\label{eqn:dba_ops_from_e-r-pairs}
   x \sqcap y &:= e(r(x) \wedge_p r(y)) &&\neg x := e(\neg_p r(x)) \\ 
   x \sqcup y &:= e'(r'(x) \vee_q r'(y)) && \lrcorner x := e'(\neg_q r'(x)). \nonumber 
\end{align}
%
%
We also define a relation $\sqsubseteq$ on $P \cup Q$ as in Equation~\ref{Eq:dBa quasiorder}. 
%
%
\noindent
We will prove that $\mathbf{P+Q}$ is a generalized glued sum (Proposition~\ref{prop:PQ-generalized glued sum}). To better understand this structure, we define a \emph{generalized D-core algebra} as a universal algebra satisfying some conditions in Definition~\ref{def:D-core+}, and show that 
\begin{itemize}
    \item  $\mathbf{P+Q}$ is a generalized D-core algebra whenever the intersection $P \cap Q$ contains both the top element of $P$ and the bottom element of $Q$.
\end{itemize}
Finally, using this result, we can characterize pure and trivial dBas. 

To find the appropriate Stone space we consider the product of two Stone spaces defined on the sets of primary ideals $\mathcal{I}_{pr}(\mathbf{D})$ and primary filters $\mathcal{F}_{pr}(\mathbf{D})$ of a dBa $\textbf{D}$. For each $x\in D$, two sets $F_x$ and $I_x$ are defined (Eq.~(\ref{eqn:FxIx})), and  
shown to be clopen subsets of the Stone space defined on  $\mathcal{F}_{pr}(\mathbf{D})$ and $\mathcal{I}_{pr}(\mathbf{D})$, respectively. 
Hence, each product $F_{x}\times I_{x}$ is a clopen subset of the product Stone space. We consider the set $\mathcal{D}:=\{(F_{x},I_{x})\mid x\in D\}$, noting that the correspondence  yields a bijection between $\mathcal{D}$ and the family of clopen rectangles of this form. Two special subfamilies are then defined: $\mathcal{D}_{\sqcap}:=\{(F_{x},I_{x\sqcap x})\mid x\in D\}$ and $\mathcal{D}_{\sqcup}:=\{(F_{x\sqcup x},I_{x})\mid x\in D\}$, each of which is shown to form a Boolean algebra under appropriate  operations. We then define two pairs of maps $r:\mathcal{D}\rightleftharpoons\mathcal{D}_{\sqcap}:e$ and $r^{\prime}:\mathcal{D}\rightleftharpoons\mathcal{D}_{\sqcup}:e^{\prime}$ and verify that both pairs satisfy the hypotheses of the New Boolean Representation Theorem (Theorem~\ref{newdbaboolean}). This yields the resulting dual Boolean algebra of clopen subsets of the product Stone space, after which we prove the following result.

   \begin{itemize}
       \item  For every dBa $\mathbf{D}$, there exists a homomorphism from $\mathbf{D}$ to the dBa of clopen subsets of a Stone space that preserves the quasi-order. If $\mathbf{D}$ is contextual, this map is an isomorphism.

   \end{itemize}

   Building on the structural properties of D-core algebras, we  introduce a sequent calculus \textbf{L} for \emph{contextual D-core algebras} and an hypersequent calculus \textbf{HL} for  \emph{pure D-core algebras}. The proof system we present is closely related to the system \textbf{CDBL} developed for contextual double Boolean algebras in \cite{howlader2021dbalogic} and \textbf{PDBL} developed for pure double Boolean algebra \cite{HOWLADER2023115} respectively, but it is simplified for clarity and ease of use. It is also proved that the systems $\textbf{L}$ and \textbf{CDBL} are equivalent while $\textbf{HL}$ and \textbf{PDBL} are equivalent.

   The paper is arranged as follows. In Section~\ref{pre}, we present the preliminaries required for our work. The notion of D-core algebras is introduced and shown to be equivalent to dBa in Section~\ref{sec:fm}. In Section~\ref{sec:re}, we present the Boolean and topological representations of dBas. In Section~\ref{see:logic}, we propose simplified proof systems \textbf{L} and \textbf{HL}  for contextual and pure D-core algebras, respectively, and show their equivalence to \textbf{CDBL} and \textbf{PDBL} as mentioned above. We conclude our work in Section~\ref{see:concl}.

\section{Preliminaries}
\label{pre}
In this section, we present  preliminaries related to dBas, including existing algebraic and topological representation results.  Our main reference are \cite{birkhoff1940lattice,wille,MR4566932,kembang2023simple,bmlpl}. We start by recalling the Boolean representation of dBas. Throughout the section, $\mathbf{ P} := (P; \wedge_p, \vee_p, \neg_{p}, \bot_p, \top_p)$  and $ 
 \mathbf{Q} := (Q; \wedge_q, \vee_q, \neg_{q},  \bot_q, \top_q)$ denote two Boolean algebras.

 Let $r:A \rightleftharpoons P: e$ and $r^{\prime}:A \rightleftharpoons Q: e^{\prime}$ be  pairs of maps such that $r\circ e=id_{P}$ and $r^{\prime}\circ e^{\prime}=id_{Q}$. We will call such pairs \emph{embedding-retraction} pairs. We define a universal algebra $\textbf{A}:=(A, \sqcap, \sqcup, \neg,\lrcorner, e^{\prime}(\top_q), e(\bot_p))$ with $\sqcap,\sqcup,\neg,\lrcorner$ defined as in Equations~\ref{eqn:dba_ops_from_e-r-pairs}.
 %
Then, we have the following result.
 
 \begin{theorem}
\cite{bmlpl} 	\label{dbaboolean}{\rm $\textbf{A}$ is a dBa iff the following holds
 		\begin{enumerate}
 		\item $e\circ r\circ e^{\prime}\circ r^{\prime}=e^{\prime}\circ r^{\prime}\circ e\circ r$
 		\item $\left\{
        \begin{aligned}
 		  e(r(x)\wedge_{p} r(e^{\prime}(r^{\prime}(x)\vee_{q} r^{\prime}(y))))&=e(r(x))  \\
          e^{\prime}(r^{\prime}(x)\vee_{q} r^{\prime}(e(r(x)\wedge_{p} r(y))))&=e^{\prime}(r^{\prime}(x))
 		\end{aligned}\right.$ for all $x, y\in A$
 		\item  $r(e^{\prime}(\top_{q}))=\top_{p}$ and $r^{\prime}(e(\bot_{p}))=\bot_{q}$.
 		\end{enumerate}
 		Moreover every dBa can be obtained from such a construction.}
 \end{theorem}
 
 \noindent Let us recall the definition of glued sum of two posets.
 \begin{definition}{\cite{birkhoff1940lattice}}\label{def:ordinal-glued-sum}
 Let \((P, \leq_P)\) and \((Q, \leq_Q)\) be two disjoint posets.
 	\begin{enumerate}
 		\item The \emph{linear sum} \(P + Q\) is the poset \((P \cup Q, \leq)\), where for \(x, y \in P \cup Q\), \(x \leq y\) if:
 		\begin{enumerate}
 			\item \(x, y \in P\) and \(x \leq_P y\),
 			\item \(x, y \in Q\) and \(x \leq_Q y\),
 			\item \(x \in P\) and \(y \in Q\).
 		\end{enumerate}
 		\item If \(P\) has a top element \(\top_p\) and \(Q\) has a bottom element \(\bot_q\), the \emph{glued sum} \(P \,\dot{+}\, Q\) is formed from \(P + Q\) by identifying \(\top_p = \bot_q\).
 	\end{enumerate}
 \end{definition}
 
 \noindent
  The glued sum of the Boolean algebras {\bf P} and {\bf Q} is 
  $(P \,\dot{+}\, Q; \sqcap, \sqcup, \neg, \lrcorner, \bot, \top) $ with  \(\bot := \bot_p\), \(\top := \top_q\) and for $x, y\in P \,\dot{+}\, Q$
 the operations are defined by:
 %
 \begin{align}\label{eqn:glued-sum Bas}
 \neg x := 
 \begin{cases}
 	\neg_p x& \text{if } x \in P, \\
 	\bot_p & \text{otherwise};
 \end{cases}
 &\quad
 \lrcorner x := 
 \begin{cases}
 	\neg_q x& \text{if } x \in Q, \\
 	\top_q & \text{otherwise.}
 \end{cases}\\
 x \sqcap y := 
 \begin{cases}  
 	x \wedge_p y & \text{if } x, y \in P, \\
 	\top_p = \bot_q & \text{if } x, y \in Q, \\
 	x & \text{if } x \in P,\, y \in Q;
 \end{cases}  &\quad 
 x \sqcup y := 
 \begin{cases}
 	x \vee_q y & \text{if } x, y \in Q, \\
 	 \bot_q = \top_p & \text{if } x, y \in P, \\
 	y & \text{if } x \in P,\, y \in Q;
 \end{cases}\nonumber
 \end{align}
 \begin{theorem}\label{thm:pure-trivial-dBa}\cite{kembang2023simple}
 	\leavevmode
 	{\rm  The glued sum $\textbf{D}:=(P \,\dot{+}\, Q, \sqcap, \sqcup, \neg, \lrcorner, \bot, \top) $ of two Boolean algebras 
 		 $\mathbf{P}$ and   $\mathbf{Q}$ is a pure and trivial dBa. Moreover every pure and trivial dBa is obtained from such construction. 
 	}
 \end{theorem}
\noindent
The next proposition lists some properties of trivial dBas to be used later.
\begin{proposition}
	\cite{kembang2023simple} \label{trivalpp}{\rm For a trivial dBa, \textbf{D} the following holds, 
	\begin{enumerate}
		\item If $x, y\in D_{\sqcap}$ then $x\sqcup y=\bot\sqcup \bot$ and $\lrcorner x=\top$.
		\item If $x, y\in D_{\sqcup}$ then $x\sqcap y=\top\sqcap \top$ and $\neg x=\bot$.
		\end{enumerate} }
\end{proposition}

Now we note the primary filter and ideal of a dBa $\textbf{D}$.

\begin{definition}
	\cite{wille, kwuida2007prime}~{\rm A   {\it filter} in  $\textbf{D}$ is a subset $F$ of $D$ such that $x\sqcap y\in F$ for all $x,y \in F$, and for all $z\in D$ and $ x \in F, x\sqsubseteq z$ implies that $z\in F$. An {\it ideal} in a dBa is defined dually.
		A filter $F$ (resp. ideal $I$) is  {\it proper} if and only if  $F\neq D$ (resp. $I\neq D$). 
		A {\it primary} filter  $F$ (resp. ideal $I$) is a non empty proper filter (resp. ideal)  such that $x\in F~ \mbox{or} ~\neg x \in F$ (resp. $~x\in I ~\mbox{or}~ \lrcorner x\in I$), for all $x\in D$.
	}
\end{definition} 

The set of primary filters is denoted by $\mathcal{F}_{pr}(\mathbf{D})$ and the set of primary ideals is denoted by $\mathcal{I}_{pr}(\mathbf{D})$. In \cite{wille}, Wille considers the following  standard context $\mathbb{K}(\mathbf{D}):=(\mathcal{F}_{pr}(\mathbf{D}), \mathcal{I}_{pr}(\mathbf{D}),\Delta)$, where for $F\in \mathcal{F}_{pr}(\mathbf{D})$ and $I\in \mathcal{I}_{pr}(\mathbf{D})$, $F\Delta I$ if and only if $F\cap I\neq\emptyset$. For any $x \in \mathbf{D}$, we define
\begin{align}\label{eqn:FxIx}
    F_{x}:=\{F\in\mathcal{F}_{pr}(\mathbf{D})\mid x\in F\} \text{ and } I_{x}:=\{I\in\mathcal{I}_{pr}(\mathbf{D})\mid x\in I\}. 
\end{align}

\begin{lemma}{\rm \cite{wille}
		\label{derivation} 
		Let $\textbf{D}$ be a dBa and  $\mathbb{K}(\textbf{D}):=(\mathcal{F}_{pr}(\textbf{D}),\mathcal{I}_{pr}(\textbf{D}),\Delta)$. 
		\begin{enumerate}
			\item $F_{x}^{\prime}=I_{x}=I_{x\sqcup x}$ for all $x\in D_{\sqcap}$.
			\item $I_{y}^{\prime}=F_{y}=F_{y\sqcap y}$ for all $y\in D_{\sqcup}$.
			\item  $F_{x}^{\prime}=I_{x\sqcap x}=I_{x_{\sqcap\sqcup}}$ and $I_{x}^{\prime}=F_{x\sqcup x}=F_{x_{\sqcup\sqcap}}$ for $x\in D$.
			\item  $(F_{x})^{c}=F_{\neg x}$ and $(I_{x})^{c}=I_{\lrcorner x}$ for $x\in D$. 
			\item $I_{x}\cap I_{y}=I_{x\sqcup y}$ and $I_{x_{\sqcup} }= I_{x}$ for $x, y\in D$.
			\item $F_{x}\cap F_{y}= F_{x\sqcap y}$ and $F_{x_{\sqcap }}=F_{x}$ for $x, y\in D$.
		\end{enumerate}
}\end{lemma}
\noindent The following representation result is then proved.
\begin{theorem}\cite{wille}
	\label{firstrepsentation}{\rm For any dBa $\textbf{D}$ and $\mathbb{K}(\textbf{D})=(\mathcal{F}_{pr}(\mathbf{D}), \mathcal{I}_{pr}(\mathbf{D}), \Delta)$, the following hold
	\begin{enumerate}
		\item For $x\in D$,  $(F_{x}, I_{x})$ is a protoconcept of $\mathbb{K}(\textbf{D})$.
		\item The map $h:\textbf{D}\rightarrow  \mathfrak{P}(\mathbb{K}(\textbf{D})) $, $x\mapsto h(x):=(F_{x}, I_{x})$ 
        is a quasi-embedding.
		\end{enumerate}}
\end{theorem} 


D\"{u}ntsch et al. \cite{duntsch2002modal} defined sufficiency, dual sufficiency, possibility and necessity operators based on a context. For a context  $\mathbb{K}:=(G,M,I)$, $g\in G$, and $m\in M$, the sets
 \begin{align*}
      g^\prime = \{m\in M\mid (g,m)\in I\} \text{ and } 
      m^\prime =\{g\in G\mid (g,m)\in I\} 
   \end{align*}
    are the {\it right-} and {\it left-neighborhood} of $g$ and $m$, respectively. For $A\subseteq G$, and $B\subseteq M$, the pairs of dual approximation operators are defined as:
\begin{align*}
   \text{possibility:} && B_{I}^{\diamondm} &:=
   \{g\in G\mid g^\prime \cap B\neq \emptyset\} \text{ and }
    A_{I}^{\diamondg} := \{m\in M\mid m^\prime \cap A\neq \emptyset\}\\
   \text{necessity:} && B_{I}^{\boxm}&:=\{g\in G\mid g^\prime \subseteq B\} \text{ and }
   A_{I}^{\boxg} := \{m\in M\mid m^\prime \subseteq A\}.
\end{align*}
%
%
   If there is no confusion about the relation involved, we shall omit 
   the subscript and denote $B_{I}^{\diamondm}$ by  $B^{\diamondm}$, $B_{I}^{\boxm}$ by $B^{\boxm}$ and  similarly for the case of $A$.

   The necessity and possibility operators correspond to approximation operators in rough set theory (RST) \cite{pawlak1982rough}. 
The following definition of object-oriented semiconcept and object-oriented protoconcept of a context $\mathbb{K}=(G, M, I)$ are adopted in \cite{howlader2018algebras,howlader2020}.

\begin{definition}
	{\rm \cite{howlader2018algebras, howlader2020} For $A\subseteq G$ and $B\subseteq M$, $(A, B)$ is an {\it object oriented semiconcept} of $\mathbb{K}$ if $A^{\boxg}=B$ or $B^{\diamondm}=A$.
		$(A, B)$ is an {\it object oriented protoconcept} of $\mathbb{K}$ if $A^{\boxg\diamondm}=B^{\diamondm}$.
}\end{definition}

$\mathfrak{R}(\mathbb{K})$  denotes the set of all object oriented protoconcepts, while the set of all object oriented semiconcepts is denoted by $\mathfrak{S}(\mathbb{K})$.
 Recall the notations $\mathfrak{P}(\mathbb{K})$ and  $\mathfrak{H}(\mathbb{K})$ for the set of protoconcepts and the set of semiconcepts of a context $\mathbb{K}$, respectively. Then we have the following results, 

 \begin{theorem} \cite{howlader2020, howlader2018algebras}
    \label{connectingthm} {\rm Let $\mathbb{K}=(G, M, I)$ be a context and let $\mathbb{K}^{c}:=(G, M, I^{c})$ be its  complemented context. i.e. $I^{c}:=G\times M\setminus I$.  Then, the following are true
     \begin{enumerate}
     \item $A^{\boxg}_{I}=A^{c\prime}_{I^{c}}$, $B^{\boxm}_{I}=B^{c\prime}_{I^{c}}$, $A^{\diamondg}_{I}=A_{I^{c}}^{\prime c}$ and $B^{\diamondm}_{I}=B_{I^{c}}^{\prime c}$.
         \item  $(A, B)\in \mathfrak{P}(\mathbb{K})$ if and only if  $(A^{c}, B)\in \mathfrak{R}(\mathbb{K}^{c})$
         \item $(A, B)\in\mathfrak{H}(\mathbb{K})$ if and only if  $(A^{c}, B)\in\mathfrak{S}(\mathbb{K}^{c})$
     \end{enumerate}}
 \end{theorem}
  
 Instead of the standard context defined by Wille, its  complement, defined by $\mathbb{K}_{pr}(\textbf{D}):=(\mathcal{F}_{pr}(\mathbf{D}), \mathcal{I}_{pr}(\mathbf{D}),\nabla \ )$ where  for $F\in \mathcal{F}_{pr}(\mathbf{D})$ and $I\in \mathcal{I}_{pr}(\mathbf{D})$, $F\nabla I$ if and only if $F\cap I=\emptyset$, is considered in \cite{MR4566932}. 
$\mathbb{K}_{pr}(\textbf{D})$ 
 is extended by equipping the sets \(\mathcal{F}_{pr}(\mathbf{D})\) 
 and \(\mathcal{I}_{pr}(\mathbf{D})\) 
 with topologies, resulting in a structure called a \emph{context on topological spaces} (CTS) \cite{MR4566932}. 
 \begin{definition}
   \cite{MR4566932}  \label{cTS} {\rm \(\mathbb{K}^{T} := ((G, \rho), (M, \tau), R)\) is called a \emph{context on topological spaces} (CTS) if \((G, M, R)\) is a formal context and \((G, \rho)\) and \((M, \tau)\) are topological spaces.}
 \end{definition}
 

In the next definition, we introduce two sets that will serve as the universes of discourse for the concrete dBas used in the representation theorem. 

\begin{definition}
\cite{MR4566932} For a CTS $\mathbb{K}^{T}$, an  object-oriented protoconcept (resp. semiconcept) $(A,B)$ of $\mathbb{K}$ is said to be  \emph{clopen } if  
$A$ is clopen in $(G, \rho)$ and $B$ is clopen in $(M, \tau)$. 
The set of all clopen object-oriented protoconcepts (resp. semiconcepts) is denoted by $\mathfrak{R}^{T}(\mathbb{K}^{T})$ (resp. $\mathfrak{S}^{T}(\mathbb{K}^{T})$).
\end{definition}

\noindent
 A topology \(\mathcal{T}\) on \(\mathcal{F}_{pr}(\mathbf{D})\) is defined by taking the family \(\mathcal{B}_0 := \{F_x : x \in D\}\) as a subbase for the closed sets, while a topology \(\mathcal{J}\) on \(\mathcal{I}_{pr}(\mathbf{D})\) is similarly generated from the subbase \(\mathcal{B} := \{I_x : x \in D\}\). It is then proved that, 
 
 \begin{proposition}\cite{MR4566932}
    \label{protopo} {\rm \((\mathcal{F}_{pr}(\mathbf{D}), \mathcal{T})\) and \((\mathcal{I}_{pr}(\mathbf{D}), \mathcal{J})\) are compact and totally disconnected topological spaces, and hence Hausdorff.}
 \end{proposition} 
 
 This result lays the topological foundation for the dual representation of pure dBas within the CTS framework. In particular, a special class of CTS, called a \emph{context on topological spaces with continuous relations} (CTSCR), is considered. For that we first recall the continuity of a relation.

\begin{definition}
\cite{MR4566932} Let $(X,\rho)$ and $(Y,\tau)$ be topological spaces. A relation $R \subseteq X \times Y$ 
is said to be \emph{continuous} if the following conditions hold:
\begin{enumerate}
   \item[(i)] If $A$ is open in $(Y,\tau)$, then both $A^{\diamondm}$ and $A^{\boxm}$ are open in $(X,\rho)$.
    \item[(ii)] If $A$ is closed in $(Y,\tau)$, then both $A^{\diamondm}$ and $A^{\boxm}$ are closed in $(X,\rho)$.
\end{enumerate}
\end{definition}

\begin{remark}
Continuity of $R^{-1}$ can equivalently be characterized using the corresponding operators induced by $R$.
\end{remark}
 Now we are ready to define a CTSCR.
 \begin{definition} \cite{MR4566932} \label{def:CTSCR topcon}
 	{\rm 
 		\(\mathbb{K}^{T} := ((G, \rho), (M, \tau), R)\) is called a \emph{context on topological spaces with continuous relations} (CTSCR) if
            the relations \(R\) and \(R^{-1}\) are continuous with respect to \((G, \rho)\) and \((M, \tau)\).
 	}
 \end{definition}
\noindent The operations $\sqcap$, $\sqcup$, $\lrcorner$, $\neg$,$\top$,$\bot$ are defined 
for $(A,B), (C,D)\in \mathfrak{R}^{T}(\mathbb{K}^{T})$ by 
\begin{center}
$(A, B)\sqcap (C, D):=(A\cup C, (A\cup C)^{\boxg})$\\
$(A, B)\sqcup (C, D):=((B\cap D)^{\diamondm}, B\cap D)$,\\
$\lrcorner(A, B):=(B^{c\diamondm} B^c)$\\ $\neg(A, B):=(A^c, A^{c\boxg})$\\
$\top :=(\emptyset, \emptyset)$~and ~$\bot :=(G, M).$
\end{center}
 
\noindent  The following results are then proved. 
 \begin{theorem}
\cite{MR4566932} 	\label{object-proto and proto}
 	{\rm  For any CTSCR $\mathbb{K}^{T} := ((G, \rho), (M, \tau), R)$, $\mathfrak{S}^{T}(\mathbb{K}^{T})\subseteq \mathfrak{R}^{T}(\mathbb{K}^{T})$. Moreover we have the following.
 		\noindent \begin{enumerate}
 			\item $\underline{\mathfrak{R}}^{T}(\mathbb{K}^{T}):=(\mathfrak{R}^{T}(\mathbb{K}^{T}),\sqcup,\sqcap,\neg,\lrcorner,\top,\bot)$ is  a fully contextual dBa.
 			\item  $\underline{\mathfrak{S}}^{T}(\mathbb{K}^{T}):=(\mathfrak{S}^{T}(\mathbb{K}^{T}),\sqcup,\sqcap,\neg,\lrcorner,\top,\bot)$ is a subalgebra of  $\underline{\mathfrak{R}}^{T}(\mathbb{K}^{T})$ and is  a pure dBa.
 			
 	\end{enumerate}}
 \end{theorem}
 
\noindent On the other hand, the following theorem can be proved.
\begin{theorem}
 	\cite{MR4566932} \label{topcontext}
 	{\rm  $\mathbb{K}_{pr}^{T}(\textbf{D}):=((\mathcal{F}_{pr}(\textbf{D}),\mathcal{T}),(\mathcal{I}_{pr}(\textbf{D}),\mathcal{J}),\nabla)$ is a CTSCR.}
 \end{theorem}

\noindent Then, we have the following representation results.
 
 \begin{theorem}[Representation theorem for dBas and contextual dBas]\cite{MR4566932} \label{thm:rep}
 \noindent{\rm\begin{enumerate}
 		\item For any dBa \(\mathbf{D}\), the map \(h : \mathbf{D} \to \underline{\mathfrak{R}}^T(\mathbb{K}_{pr}^T(\mathbf{D}))\) defined by $h(x) := (F_{\neg x}, I_x)$ for all  $x \in \mathbf{D}$,
 		is a quasi-embedding from \(\mathbf{D}\) into \(\underline{\mathfrak{R}}^T(\mathbb{K}_{pr}^T(\mathbf{D}))\). Moreover, \(\mathbf{D}_p\) is isomorphic to \(\underline{\mathfrak{S}}^T(\mathbb{K}_{pr}^T(\mathbf{D}))\).
 		
 		\item For any contextual dBa \(\mathbf{D}\), the above map \(h\) is an embedding from \(\mathbf{D}\) into \(\underline{\mathfrak{R}}^T(\mathbb{K}_{pr}^T(\mathbf{D}))\).
 	\end{enumerate}}
 \end{theorem}

 \begin{theorem}[\textbf{Representation theorem for fully contextual dBa}]\cite{MR4566932}
 	\label{iso-fullycxt dBa}
 	{\rm Any fully contextual dBa $\textbf{D}$ is isomorphic to $\underline{\mathfrak{R}}^{T}(\mathbb{K}_{pr}^{T}(\textbf{D}))$, the algebra of clopen object-oriented protoconcepts.}
 \end{theorem}
 
  \begin{theorem}[\textbf{Representation theorem for pure dBas}]\cite{MR4566932}
 	\label{RTDBA} 
 	{\rm Any pure dBa $\textbf{D}$ is isomorphic to  $\underline{\mathfrak{S}}^{T}(\mathbb{K}_{pr}^{T}(\textbf{D}))$, the algebra of clopen  object-oriented semiconcepts.}
 \end{theorem}



\section{Elimination of Redundant Axioms}\label{sec:fm}
 In this section, we extract a \emph{minimal set of axioms} equivalent to the set of equations in Definition~\ref{def:DBA}. In particular, we proposed the following definition.
\begin{definition} \label{def:D-core+}
    {\rm  An  algebra $ \textbf{D}:= (D; \sqcup, \sqcap, \neg,\lrcorner,\top,\bot)$ satisfying the following properties is called a  {\it D-core algebra}. For any $x,y,z \in D$,\\
    $\begin{array}{ll}
			(1a) \  x\sqcap y = y\sqcap  x  &
			(1b) \   x \sqcup   y = y\sqcup   x  \\
			(2a) \  \neg (x \sqcap  x) = \neg  x  &
			(2b) \   \lrcorner(x \sqcup   x )= \lrcorner x \\
			(3a) \   x  \sqcap (x \sqcup y)=x \sqcap  x  &
			(3b) \   x \sqcup  (x \sqcap y) = x \sqcup   x \\
			(4a) \  x \sqcap  (y \vee z ) = (x\sqcap  y)\vee (x \sqcap  z) &
			(4b) \   x \sqcup  (y \wedge z) = (x \sqcup  y) \wedge  (x \sqcup  z) \\
			(5a) \   \neg \neg (x \sqcap  y)= x \sqcap  y &
			(5b) \   \lrcorner\lrcorner(x \sqcup  y) = x\sqcup  y \\
			(6a) \   x  \sqcap \neg  x= \bot &
			(6b) \   x \sqcup \lrcorner x = \top  \\
			(7) \   (x \sqcap  x) \sqcup (x \sqcap x) = (x \sqcup x) \sqcap (x \sqcup x) &
		\end{array}$}
\end{definition}

In the rest of this section, we will prove that Definition \ref{def:D-core+} is equivalent to 
Definition \ref{def:DBA}. In the following results, the equations come in pairs, with (b) being the dual of (a).
Therefore, we will in the sequel only write the proof for (a)'s and get (b)'s as dual. 

\begin{proposition}\label{prop:axiom1a}
    {\rm Let $\textbf{D}$ be a D-core algebra and  $x, y, z\in D$. Then, the following hold.

    \noindent$\begin{array}{ll}
    (1a)~ x\sqcap x=\neg\neg x & (1b)~ x\sqcup x=\lrcorner\lrcorner x\\
    (2a)~ (x\sqcap y)\sqcap (x\sqcap y)=x\sqcap y & (2b)~ (x\sqcup y)\sqcup (x\sqcup y)=x\sqcup y\\
    (3a)~ \neg (x\sqcap (y\vee  z))=\neg (x\sqcap y)\sqcap \neg (x\sqcap z) & (3b)~ \lrcorner (x\sqcup ( y\wedge  z))=\lrcorner (x\sqcup y)\sqcup \lrcorner (x\sqcup z)\\
    (4a)~ x\vee x= x\sqcap x & (4b)~ x\wedge x=x\sqcup x \\
    (5a)~ \neg\neg\neg x = \neg x & (5b)~ \lrcorner\lrcorner\lrcorner x = \lrcorner x \\
     (6a)~\neg x\sqcap \neg x = \neg x & (6b)~\lrcorner x\sqcup \lrcorner x = \lrcorner x\\
    \end{array}$ 

    }
\end{proposition}
\begin{proof} Let $x,y,z\in D$.

\noindent $(1a)$  Setting $y:=x$ in Def.~\ref{def:D-core+}~(5a) and then using Def.~\ref{def:D-core+}~(2a) we get \footnote{We refer to definitions, propositions, lemmas or theorems with the short form Def., Prop., Lem. or Thm., (or even D.,P.,T. or L. inside equations) followed by the identifying number and the relevant numbering of the identity, if applies. The commutativity of $\sqcap,\sqcup$ (Def.~\ref{def:D-core+}~(1a),(1b)) and of $\wedge,\vee$, will be used without explicitly referenced.}
\begin{align*}
 x\sqcap x \overset{\mathrm{D.\ref{def:D-core+}}~(5a)}{=} \neg\neg(x\sqcap x) \overset{\mathrm{D.\ref{def:D-core+}}~(2a)}{=} \neg\neg x, \text{ and (1a) is proved.}    
\end{align*}
%
  \noindent  $(2a)$ Applying Prop.~\ref{prop:axiom1a}~(1a) to Def.~\ref{def:D-core+}~(5a) we get 
  \begin{align*}
  x\sqcap y \overset{{\rm D.}\ref{def:D-core+}~(5a)}{=} \neg\neg(x\sqcap y) \overset{{\rm P.}\ref{prop:axiom1a}~(1a)}{=} (x\sqcap y)\sqcap (x\sqcap y), \text{ and (2a) is proved.}
  \end{align*}
  %
%
    \noindent $(3a)$ $x\sqcap (y\vee  z) = \neg (\neg (x\sqcap y)\sqcap \neg(x\sqcap z))$ by Def.~\ref{def:D-core+}~(4a) and the definition of $\vee$. Taking negation $\neg$ on both sides and using Def.\ref{def:D-core+}~(5a) we get
    \begin{align*}
 \neg( x\sqcap (y\vee z)) = \neg\neg (\neg (x\sqcap y)\sqcap \neg(x\sqcap z)) \overset{(5a)}{=}
      \neg (x\sqcap y)\sqcap \neg(x\sqcap z).
    \end{align*}
    

    \noindent $(4a)$ $x\vee x \overset{\mathrm{def.}\vee}{=} \neg (\neg x \sqcap \neg x) \overset{{\rm D.}\ref{def:D-core+}~(2a)}{=} \neg \neg x \overset{{\rm P.}\ref{prop:axiom1a}~(1a)}{=} x\sqcap x.$
         


     \noindent 
     (5a) $\begin{aligned}
         \neg x \overset{{\rm D.}\ref{def:D-core+}~(2a)}{=} \neg(x\sqcap x) \overset{{\rm D.}\ref{def:D-core+}~(5a)}{=} \neg(\neg\neg(x\sqcap x)) \overset{{\rm D.}\ref{def:D-core+}~(2a)}{=} \neg\neg\neg x.
     \end{aligned}
     $
     
\noindent (6a)\quad 
$\begin{aligned}
    \neg x \sqcap \neg x \overset{{\rm P}\ref{prop:axiom1a}~(1a)}{=}\neg\neg\neg x \overset{{\rm P}\ref{prop:axiom1a}~(5a)}{=} \neg x. 
\end{aligned}$

\end{proof}
\noindent
In Proposition~\ref{prop:axiom1a}, the identities (1a) and (1b) can be interpreted as double negation laws, and (2a), (2b) as idempotency in $D_\sqcap$ and $D_\sqcup$.  
\noindent The following theorem demonstrates that the axioms $(1a),(1b),  (9a), (9b), (11a)$ and $(11b)$ of double Boolean algebra (Def.~\ref{def:DBA}) are derivable within D-core~algebras.
\begin{theorem} \label{thm:firstindeax}
	{\rm   Let $\textbf{D}$ be a D-core algebra and  $x, y, z\in D$. Then, the following hold.

   \noindent  $\begin{array}{ll}
          (1a)~ (x \sqcap x ) \sqcap  y = x \sqcap  y & (1b)~ (x \sqcup x)\sqcup  y = x \sqcup y\\
         (2a)~ \neg \top=\bot & (2b)~ \lrcorner \bot=\top\\
         (3a)~ \neg \bot = \top \sqcap   \top & (3b)~ \lrcorner\top =\bot \sqcup  \bot
         \end{array}$
		}
\end{theorem}

\begin{proof} 
Let $x, y, z\in D$.


\noindent	(1a) \quad 
${\begin{aligned}[t]
(x\sqcap x)\sqcap y &\overset{{\rm D.}\ref{def:D-core+}~(1a)}{=} y\sqcap(x\sqcap x) \overset{{\rm P.}\ref{prop:axiom1a}~(4a)}{=} y\sqcap(x\vee x)  
\overset{{\rm D.}\ref{def:D-core+}~(4a)}{=} (y\sqcap x )\vee (y\sqcap x)\\
 &\overset{{\rm P.}\ref{prop:axiom1a}~(4a)}{=} (y\sqcap x)\sqcap (y\sqcap x)\overset{{\rm P.}\ref{prop:axiom1a}~(2a)}{=} y\sqcap x
 \overset{{\rm D.}\ref{def:D-core+}~(1a)}{=} x\sqcap y.
\end{aligned}}$

\noindent (2a) \quad ${\begin{aligned}[t]
\neg \top & 
\overset{{\rm P}\ref{prop:axiom1a}~(5a)}{=} \neg\neg\neg \top   
\overset{{\rm P}\ref{prop:axiom1a}~(1a)}{=} \neg\top\sqcap\neg\top 
\overset{{\rm D}\ref{def:D-core+}~(3a)}{=} \neg\top\sqcap(\neg\top \sqcup \lrcorner\neg\top)\\ 
&\overset{{\rm D}\ref{def:D-core+}~(6b)}{=} \neg\top \sqcap \top 
\overset{{\rm D}\ref{def:D-core+}~(6a)}{=} \bot.
\end{aligned}}$
 	
\noindent (3a)\quad ${\begin{aligned}
	     \top\sqcap \top &\overset{{\rm P}\ref{prop:axiom1a}~(1a)}{=} \neg\neg \top &\overset{{\rm T}\ref{thm:firstindeax}~(2a)}{=} \neg \bot.
	 \end{aligned}}$ \end{proof}
The next proposition provides the De Morgan laws for a D-core algebra.
\begin{proposition}
   \label{prop: De Morgan}{\rm Let $\textbf{D}$ be a D-core algebra and  $x, y \in D$. Then, it holds:\\
    \noindent  $\begin{array}{ll}
         (1a)~\neg(x\sqcap y)=\neg x\vee\neg y& (1b)~\lrcorner(x\sqcup y)=\lrcorner x\wedge\lrcorner y\\
         (2a)~\neg(x\vee y)=\neg x\sqcap \neg y&(2b)~\lrcorner(x\wedge y)=\lrcorner x\sqcup \lrcorner y
    \end{array}$}
\end{proposition}
\begin{proof}
     Let $x, y \in D$
     
\noindent (1a)\quad 
$\begin{aligned}[t]
&\neg(x\sqcap y)
   \overset{\text{T.\,\ref{thm:firstindeax}\,(1a)}}{=}
   \neg((x\sqcap x)\sqcap (y\sqcap y))
\overset{\text{P.\,\ref{prop:axiom1a}\,(1a)}}{=}
   \neg(\neg\neg x\sqcap \neg\neg y)
\overset{\text{def.}\vee}{=}
   \neg x \vee \neg y
\end{aligned}$

\noindent(2a)\quad 
$\begin{aligned}[t]
&\neg(x\vee y)
   \overset{\text{def.}\vee}{=}
   \neg\neg(\neg x \sqcap \neg y)\overset{{\rm D.}\ref{def:D-core+}\,(5a)}{=}
   \neg x \sqcap \neg y
\end{aligned}$

\end{proof}
\noindent
The next proposition computes some meet and join with top and bottom. 
\begin{proposition} \label{prop:axion6}
{\rm   Let $\textbf{D}$ be a D-core algebra and  $x, y\in D$. We have:

   \noindent  $\begin{array}{ll}
    (1a)~x\sqcap \top=x\sqcap x & (1b)~x\sqcup \bot=x\sqcup x\\
    (2a)~\neg x \sqcap \top =\neg x & (2b)~\lrcorner x \sqcup \bot =\lrcorner x .\\
    (3a)~(x\sqcap y) \sqcap \top =  x \sqcap y & (3b)~ (x \sqcup y) \sqcap \top =  x\sqcup y\\
    (4a)~\neg x\sqcap \neg\bot=\neg x& (4b)~\lrcorner x\sqcup \lrcorner\top=\lrcorner x.\\
    (5a)~x\vee\bot = x\sqcap x & (5b)~x\wedge \top = x\sqcup x\\
    (6a)~ \neg x\vee \bot=\neg x &(6b)~ \lrcorner x\wedge \top =\lrcorner x\\
    (7a)~(x\sqcap y) \vee \bot =  x \sqcap y & (7b)~ (x \sqcup y) \wedge \top =  x\sqcup y\\
    (8a)~x\sqcap ( y\vee\top)=x\sqcap ( x\vee y) & (8b)~x\sqcup  ( y\wedge \bot)=x\sqcup (x\wedge y)\\
    (9a)~\bot\sqcap \bot=\bot & (9b)~\top\sqcup \top=\top\\
    (10a)~(\bot\sqcap \neg x)\sqcap x=\bot & (10b)~(\top\sqcup \lrcorner x)\sqcup x=\top\\ 
    (11a)~(\bot\sqcap x)\sqcap\neg x=\bot & (11b)~(\top\sqcup x)\sqcup\lrcorner x=\top\\ 
    (12a)~\bot\sqcap\lrcorner x=\bot & (12b)~\top\sqcup\neg x=\top\\
    (13a)~\bot\sqcap \neg\lrcorner x=\bot & (13b)\top\sqcup \lrcorner\neg x=\top.
    \end{array}$}
\end{proposition}
\begin{proof} For any $x, y\in D$,

\noindent	(1a) \quad 
$\begin{aligned}[t]
x\sqcap \top \overset{{\rm D.}\ref{def:D-core+}~(6b)}{=} x\sqcap (x\sqcup \lrcorner x) \overset{{\rm D.}\ref{def:D-core+}~(3a)}{=} x\sqcap x.
\end{aligned}$

\noindent (2a)\quad 
$\begin{aligned}
   \neg x \sqcap \top \overset{{\rm P}\ref{prop:axion6}~(1a)}{=} \neg x \sqcap \neg x \overset{{\rm P}\ref{prop:axiom1a}~(6a)}{=}\neg x. 
\end{aligned}
$

\noindent (3a)\quad 
$\begin{aligned}
(x\sqcap y) \sqcap \top \overset{{\rm P.}\ref{prop:axion6}~(1a)}{=} (x\sqcap y) \sqcap (x \sqcap y)  \overset{{\rm P.}\ref{prop:axiom1a}~(2a)}{=} x \sqcap y. 
\end{aligned}
$

\noindent  (4a)
$\begin{aligned}[t]
\neg x \sqcap \neg \bot
&\overset{\text{T.\,\ref{thm:firstindeax}(3a)}}{=}
\neg x \sqcap (\top \sqcap \top)
\overset{\text{T.\,\ref{thm:firstindeax}(1a)}}{=}
\neg x \sqcap \top \overset{{\rm P}\ref{prop:axion6}~(2a)}{=}
\neg x.
\end{aligned}$

\noindent(5a)\quad
${\begin{aligned}[t]
x\vee \bot
&\overset{\rm{def.}\vee}{=} \neg(\neg x \sqcap \neg \bot)
\overset{{\rm P}\ref{prop:axion6}~(4a)}{=}
\neg \neg x
\overset{{\rm P.}\ref{prop:axiom1a}~(1a)}{=}
x\sqcap x.
\end{aligned}}$

\noindent(6a)\quad
${\begin{aligned}[t]
\neg x\vee \bot &\overset{{\rm P.}\ref{prop:axion6}~(5a)}{=} \neg x \sqcap \neg x \overset{{\rm P}\ref{prop:axiom1a}~(6a)}{=}\neg x. 
\end{aligned}}$

\noindent(7a)\quad
${\begin{aligned}[t]
 (x\sqcap y) \vee \bot &\overset{{\rm P.}\ref{prop:axion6}~(5a)}{=}  (x\sqcap y) \sqcap (x\sqcap y) \overset{{\rm P}\ref{prop:axiom1a}~(6a)}{=} x\sqcap y. 
\end{aligned}}$

\noindent (8a)\quad
${\begin{aligned}[t]
	x\sqcap ( y\vee \top)&\overset{{\rm D.}\ref{def:D-core+}~(4a)}{=}(( x\sqcap  y)\vee( x\sqcap \top)) 
    \overset{{\rm P.}\ref{prop:axion6}~(1a)}{=}(( x\sqcap  y)\vee( x\sqcap x))\\
	&\overset{{\rm D.}\ref{def:D-core+}~(4a)}{=}x\sqcap ( y\vee x).
\end{aligned}}$

\noindent (9a)\quad 
$\begin{aligned}
   \bot \sqcap \bot \overset{{\rm T.}\ref{thm:firstindeax}~(2a)}{=} \neg \top \sqcap \neg \top \overset{{\rm P.}\ref{prop:axiom1a}~(6a)}{=} \neg \top \overset{{\rm T.}\ref{thm:firstindeax}~(2a)}{=}\bot . 
\end{aligned}
$


\noindent (10a)\quad
${\begin{aligned}[t]
\bot
&\overset{{\rm D.}\text{\ref{def:D-core+} (6a)}}{=}
(\bot\sqcap\neg x)\sqcap \neg(\bot\sqcap\neg x)
\overset{{\rm T.}\text{\ref{thm:firstindeax} (2a)}}{=}
(\bot\sqcap\neg x)\sqcap \neg(\neg\top\sqcap\neg x)\\[4pt]
&\overset{\text{def.}\vee}{=}
(\bot\sqcap\neg x)\sqcap (\top\vee x)
\overset{{\rm P.}\text{\ref{prop:axion6} (8a)}}{=}
(\bot\sqcap\neg x)\sqcap \big((\bot\sqcap\neg x)\vee x\big)\\[2pt]
&\overset{\text{def.}\vee}{=}
(\bot\sqcap\neg x)\sqcap \neg\big(\neg(\bot\sqcap\neg x)\sqcap \neg x\big)\\[2pt]
&\overset{{\rm T.}\text{\ref{thm:firstindeax} (2a), def.}\vee}{=}(\bot\sqcap\neg x)\sqcap \neg\big((\top\vee x)\sqcap \neg x\big)\\[2pt]
&\overset{{\rm P.}\text{\ref{prop:axion6} (8a)}}{=}(\bot\sqcap\neg x)\sqcap \neg\big(\neg x\sqcap(x\vee\neg x)\big)\\[2pt]
&\overset{{\rm D.}\text{\ref{def:D-core+}(4a)}}{=}
(\bot\sqcap\neg x)\sqcap \neg\big((\neg x\sqcap x)\vee (\neg x\sqcap\neg x)\big)\\[2pt]
&\overset{{\rm D.}\text{\ref{def:D-core+}(6a),\ref{prop:axiom1a}~(6a)}}{=}
(\bot\sqcap\neg x)\sqcap \neg(\bot\vee \neg x)
\overset{{\rm P.}\text{\ref{prop:axion6}(6a)}}{=}
(\bot\sqcap\neg x)\sqcap \neg\neg x
\\[2pt]
&\overset{{\rm P.}\text{\ref{prop:axiom1a} (1a)}}{=}
(\bot\sqcap\neg x)\sqcap (x\sqcap x)
\overset{{\rm T.}\text{\ref{thm:firstindeax} (1a)}}{=}
(\bot\sqcap\neg x)\sqcap x.
\end{aligned}}$

\noindent (11a)\quad
Setting $x:=\neg x$ in Proposition~\ref{prop:axion6}~(10a) we get

${\begin{aligned}[t]
\bot = 
(\bot\sqcap\neg\neg x)\sqcap \neg x
\overset{{\rm P.}\text{\ref{prop:axiom1a} (1a)}}{=}
(\bot\sqcap (x\sqcap x))\sqcap \neg x\overset{{\rm T.}\text{\ref{thm:firstindeax} (1a)}}{=}
(\bot\sqcap x)\sqcap \neg x.
\end{aligned}}$

 \noindent(12a)\quad
  {$\begin{aligned}
      \bot &\overset{{\rm P.}\ref{prop:axion6}~(9a)}{=}\bot\sqcap \bot\overset{{\rm D.}\ref{def:D-core+}~(3a)}{=}\bot\sqcap(\bot\sqcup \lrcorner x)\overset{{\rm P.}\ref{prop:axion6}~(2b)}{=}\bot\sqcap \lrcorner x
  \end{aligned}$}

 \noindent (13a)\quad  $\begin{aligned}[t]
  	\bot&\overset{{\rm P.}\ref{prop:axion6}~(11a)}{=}(\bot\sqcap \lrcorner x)\sqcap \neg\lrcorner x\overset{{\rm P.} \ref{prop:axion6}~(12a)}{=} \bot\sqcap \neg\lrcorner x.
  \end{aligned}$

\end{proof}

Although $\bot\sqcap x=\bot$ for $x \in \{\bot, \lrcorner x, \neg\lrcorner x\}$, it is more elaborate to prove that this holds for any $x\in D$. This is the goal of the next two propositions. 

\begin{proposition}
	\label{prop:foraxiom6}{\rm Let $\textbf{D}$ be a D-core algebra and  $x, y\in D$. The following hold:
    
	 \noindent  $\begin{array}{ll}
	(1a)~x\sqcup (\bot\sqcup y)=x\sqcup y& (1b)~x\sqcap (\top\sqcap y)=x\sqcap y\\
	(2a)~(\bot\sqcap x)\sqcup y=\bot\sqcup y & (2b)~(\top\sqcup x)\sqcap y=\top\sqcap y\\
	(3a)~(\bot\sqcap x)\sqcap (\bot\sqcup y)=\bot\sqcap x& (3b)~(\top\sqcup x)\sqcup (\top\sqcap y)=\top\sqcup x\\
	(4a)~(\bot \sqcap x)\sqcap\neg(x\sqcap \neg y)=(\bot\sqcap x)\sqcap y & (4b)~(\top \sqcup x)\sqcup\lrcorner(x\sqcup \lrcorner y)=(\top\sqcup x)\sqcup y.
	\end{array}$}
\end{proposition}
\begin{proof}
Let $x, y\in D$,

\noindent (1a)\quad$\begin{aligned}[t]
x\sqcup(\bot\sqcup y)
&\overset{\text{P.\ref{prop:axion6}(1b)}}{=}
x\sqcup (y\sqcup y)
\overset{\text{T.\ref{thm:firstindeax}(1b)}}{=}
x\sqcup y.
\end{aligned}$
    
	\noindent (2a)\quad 
$\begin{aligned}[t]
y\sqcup(\bot\sqcap x)
&\overset{\text{P.\ref{prop:foraxiom6}(1a)}}{=} y\sqcup(\bot\sqcup(\bot\sqcap x))
\overset{\text{D.\ref{def:D-core+}(3b)}}{=}
y\sqcup(\bot\sqcup\bot)
\overset{\text{T.\ref{thm:firstindeax}(1b)}}{=}
y\sqcup\bot.
\end{aligned}$
    
	\noindent (3a)\quad 
$\begin{aligned}[t]
&\phantom{\overset{{\rm D.\ref{def:D-core+} (000)}}{=}} (\bot\sqcap x)\sqcap(\bot\sqcup y)
\overset{\text{T.\,\ref{thm:firstindeax}(1b)}}{=}
(\bot\sqcap x)\sqcap\big((\bot\sqcup\bot)\sqcup y\big) \\
&\overset{\text{D.\,\ref{def:D-core+}(3b)}}{=}
(\bot\sqcap x)\sqcap\big((\bot\sqcup(\bot\sqcap x))\sqcup y\big) 
\overset{\text{P.\,\ref{prop:foraxiom6}(1a)}}{=}
(\bot\sqcap x)\sqcap\big(y\sqcup(\bot\sqcap x)\big) \\
&\overset{\text{D.\,\ref{def:D-core+}(3a)}}{=}
(\bot\sqcap x)\sqcap(\bot\sqcap x) 
\overset{\text{P.\,\ref{prop:axiom1a}(2a)}}{=}
\bot\sqcap x.
\end{aligned}$

\noindent (4a)\quad
$\begin{aligned}[t]
&\phantom{\overset{{\rm D.\ref{def:D-core+} (000)}}{=}} (\bot\sqcap x)\sqcap\neg(x\sqcap \neg y)
\overset{{\rm P.}\ref{prop: De Morgan}~(1a)}{=} (\bot\sqcap x)\sqcap (\neg x \vee \neg\neg y)\\
&\overset{{\rm D.}\ref{def:D-core+}~(4a)}{=} ((\bot\sqcap x)\sqcap \neg x) \vee ((\bot\sqcap x)\sqcap  \neg\neg y) \overset{{\rm P.}\ref{prop:axion6}~(11a)}{=} \bot \vee ((\bot\sqcap x)\sqcap  \neg\neg y) \\
& \overset{{\rm P.}\ref{prop:axion6}~(7a)}{=}(\bot\sqcap x)\sqcap  \neg\neg y 
\overset{{\rm P.}\ref{prop:axiom1a}~(1a)}{=}(\bot\sqcap x)\sqcap (y\sqcap y) 
\overset{{\rm T.}\ref{thm:firstindeax}(1a)}{=}
(\bot\sqcap x)\sqcap y \\
\end{aligned}$

\end{proof}

\begin{proposition}
	\label{proaxiom6}
	{\rm Let $\textbf{D}$ be a D-core algebra and  $x, y\in D$. The following hold:
    
	 \noindent  $\begin{array}{ll}
		(1a)~\bot\sqcap \neg(\neg x\sqcap \lrcorner y)=\bot\sqcap x &(1b)~\top\sqcup \lrcorner(\lrcorner x\sqcup \neg y)=\top\sqcup x\\
		 (2a)~x\sqcap\neg(\neg x\sqcap \neg y)=x\sqcap \neg (\neg y\sqcap \neg \top) & (2b)~x\sqcup\lrcorner(\lrcorner x\sqcup \lrcorner y)=x\sqcup \lrcorner (\lrcorner y\sqcup \lrcorner \bot)\\
		(3a)~x\sqcap \neg(\neg y\sqcap \neg(x\sqcup z))=x\sqcap \neg(\bot\sqcap \neg y)& (3b)~x\sqcup \lrcorner(\lrcorner y\sqcup \lrcorner(x\sqcap z))=x\sqcup \lrcorner(\top\sqcup \lrcorner y)\\
		(4a)~\bot\sqcap x=\bot & (4b)~\top\sqcup x=\top.
		\end{array}$}
\end{proposition}
\begin{proof} For $x, y\in D$

\noindent (1a)\quad 
$\begin{aligned}[t]
\bot\sqcap \neg(\neg x\sqcap \lrcorner y)
&\overset{\text{T.\,\ref{thm:firstindeax}(1a)}}{=}
\bot\sqcap \neg(\neg x\sqcap (\lrcorner y\sqcap \lrcorner y))
\overset{\text{P.\,\ref{prop:axiom1a}(1a)}}{=}
\bot\sqcap \neg(\neg x\sqcap \neg\neg \lrcorner y)
\\
&\overset{{\rm def.}\vee}{=} \bot \sqcap (x\vee \neg\lrcorner y) \overset{{\rm D.}\ref{def:D-core+}~(4a)}{=} (\bot \sqcap x)\vee (\bot \sqcap \neg\lrcorner y) \\
&\overset{\text{def.}\vee}{=}
\neg(\neg(\bot\sqcap x)\sqcap \neg(\bot\sqcap \neg\lrcorner y))
\overset{\text{P.\,\ref{prop:axion6}~(13a)}}{=}
\neg(\neg(\bot\sqcap x)\sqcap \neg \bot) \\
&\overset{\text{T.\,\ref{thm:firstindeax}(3a)}}{=}
\neg(\neg(\bot\sqcap x)\sqcap (\top\sqcap \top))
\overset{\text{T.\,\ref{thm:firstindeax}(1a)}}{=}
\neg(\neg(\bot\sqcap x)\sqcap \top) \\
&\overset{{\rm P.}\ref{prop:axion6}(2a)}{=}
\neg\neg(\bot\sqcap x)
\overset{\text{D.\,\ref{def:D-core+}~(5a)}}{=}
\bot\sqcap x.
\end{aligned}$



\noindent (2a)\quad 
$\begin{aligned}[t]
x\sqcap\neg(\neg y\sqcap \neg\top)
\overset{{\rm def.}\vee}{=} x \sqcap (y\vee \top) 
\overset{{\rm P.}\ref{prop:axion6}(8a)}{=} x \sqcap (y\vee x) 
\overset{{\rm def.}\vee}{=}
x\sqcap \neg(\neg x\sqcap \neg y)
\end{aligned}$


\noindent (3a) \quad
$\begin{aligned}[t]
&\phantom{\overset{{\rm D.\ref{def:D-core+} (000)}}{=}} x\sqcap \neg (\neg y\sqcap \neg(x\sqcup z))
\overset{{\rm def.}\vee }{=} x \sqcap (y\vee (x\sqcup z)) \\
&\overset{{\rm D.}\ref{def:D-core+}(4a)}{=} (x \sqcap y) \vee (x\sqcap (x\sqcup z))
\overset{{\rm D.}\ref{def:D-core+}(3a)}{=} (x \sqcap y) \vee (x\sqcap x)
\overset{{\rm D.}\ref{def:D-core+}(4a)}{=} x\sqcap (y\vee x) \\
&\overset{{\rm P.}\ref{prop:axion6}~(8a)}{=} x\sqcap (y\vee \top) 
\overset{{\rm def.}\vee}{=} 
x\sqcap \neg(\neg\top\sqcap \neg y)
\overset{\text{T.\,\ref{thm:firstindeax}(2a)}}{=}
x\sqcap \neg(\bot\sqcap \neg y)
\end{aligned}$

			
\noindent (4a)\quad 
$\begin{aligned}[t]
&\bot
\overset{\text{D.\,\ref{def:D-core+}(6a)}}{=}
(\bot\sqcap x)\sqcap \neg (\bot\sqcap x)
\overset{\text{T.\,\ref{thm:firstindeax}(1a)}}{=}
(\bot\sqcap x)\sqcap \neg (\bot\sqcap (x\sqcap x))
\\
&\overset{\text{P.\,\ref{prop:axiom1a}(1a)}}{=}
(\bot\sqcap x)\sqcap \neg (\bot\sqcap \neg\neg x)
\overset{\text{P.\,\ref{proaxiom6}(3a)}}{=}
(\bot\sqcap x)\sqcap \neg (\neg\neg x \sqcap \neg((\bot\sqcap x)\sqcup y))
\\
&\overset{\text{P.\,\ref{prop:foraxiom6}(2a)}}{=}
(\bot\sqcap x)\sqcap \neg (\neg\neg x \sqcap \neg(\bot\sqcup y))
\overset{\text{P.\,\ref{prop:axiom1a}(1a)}}{=}
(\bot\sqcap x)\sqcap \neg ((x\sqcap x)\sqcap \neg(\bot\sqcup y))
\\
&\overset{\text{T.\,\ref{thm:firstindeax}(1a)}}{=}
(\bot\sqcap x)\sqcap \neg (x\sqcap \neg(\bot\sqcup y))
\overset{\text{P.\,\ref{prop:foraxiom6}(4a)}}{=}
(\bot\sqcap x)\sqcap (\bot\sqcup y)
\overset{\text{P.\,\ref{prop:foraxiom6}(3a)}}{=}
\bot\sqcap x
\end{aligned}$
\end{proof}
\noindent
The next theorem establishes that the axioms $(6a)$ and $(6b)$ of double Boolean algebra (Def.~\ref{def:DBA}) are derivable within the framework of D-core algebra.
\begin{theorem}
\label{thm:axiom6ad}	{\rm Let $\textbf{D}$ be a D-core algebra and  $x, y, \in D$. Then, we have
\[\begin{array}{cc}
    (1a)~ x\sqcap (x\vee y)=x\sqcap x & (1b)~ x\sqcup (x\wedge y)=x\sqcup x.  
\end{array}\] }
\end{theorem}

\begin{proof}
For $x, y\in D$ 

\noindent $\begin{aligned}[t]
x\sqcap \neg(\neg x\sqcap \neg y)
&\overset{\text{P.\,\ref{proaxiom6}(2a)}}{=}
x\sqcap \neg(\neg\top\sqcap \neg y)
\overset{\text{T.\,\ref{thm:firstindeax}(2a)}}{=}
x\sqcap \neg(\bot\sqcap \neg y)\\
&\overset{\text{P.\,\ref{proaxiom6}(3a)}}{=}
x\sqcap \neg(\neg y\sqcap \neg(x\sqcup \top))
\overset{\text{P.\,\ref{proaxiom6}(4b)}}{=}
x\sqcap \neg(\neg y\sqcap \neg\top)\\
&\overset{\text{T.\,\ref{thm:firstindeax}(2a)}}{=}
x\sqcap \neg(\neg y\sqcap \bot)
\overset{\text{P.\,\ref{proaxiom6}(4a)}}{=}
x\sqcap \neg\bot\\
&\overset{\text{T.\,\ref{thm:firstindeax}(3a)}}{=}
x\sqcap (\top\sqcap \top)
\overset{\text{T.\,\ref{thm:firstindeax}(1a)}}{=}
x\sqcap \top
\overset{{\rm P.} \ref{prop:axion6}~(1a)}{=} x\sqcap x
\end{aligned}$

	\noindent
\end{proof}

\begin{proposition}
\label{prop:meet-neg}
    {\rm Let $\textbf{D}$ be a D-core algebra and  $x, y,z \in D$. Then we have:

	 \noindent  $\begin{array}{ll}
    (1a)~x\sqcap \neg (x\sqcap y)=x\sqcap \neg y& (1b)~x\sqcup \lrcorner (x\sqcup y)=x\sqcup \lrcorner y\\
     (2a)~x\sqcap\neg(x \sqcap \neg y)=x\sqcap y & (2b)~x\sqcup\lrcorner(\lrcorner y\sqcup x)=x\sqcup y\\
    (3a)~ x\sqcap \neg(\neg x\sqcap y) =  x\sqcap x &(3b)~\lrcorner((x\sqcup x)\sqcup \lrcorner(\neg x\sqcup y))=\lrcorner x\\
    (4a)~\neg x\sqcap \neg(x\sqcap y))=\neg x &(4b)~ \lrcorner x\sqcup \lrcorner(x\sqcup y))= \lrcorner x \\
    (5a)~\neg x\sqcap \neg((x\sqcap y)\sqcap z)=\neg x& (5b)~\lrcorner x\sqcup \lrcorner((x\sqcup y)\sqcup z)=\lrcorner x\\
    (6a)~x\sqcap (y\sqcap \neg x)=\bot& (6b)~ x\sqcup (y\sqcup \lrcorner x)=\top
    \end{array}$}
\end{proposition}
\begin{proof} Let $x, y,z \in D$

\noindent (1a)\quad 
$\begin{aligned}[t]
x\sqcap \neg (x\sqcap y) &\overset{{\rm P.}\ref{prop: De Morgan}~(1a)}{=}
x\sqcap (\neg x \vee \neg y) 
\overset{{\rm D.}\ref{def:D-core+}~(4a)}{=} (x\sqcap \neg x) \vee (x\sqcap \neg y) \\
& \overset{{\rm D.}\ref{def:D-core+}~(6a)}{=} \bot \vee (x\sqcap \neg y) \overset{{\rm P.}\ref{prop:axion6}~(7a)}{=} x\sqcap \neg y. 
\end{aligned}$

\noindent (2a)\quad 
$\begin{aligned}[t]
x\sqcap \neg (x\sqcap \neg y) &\overset{{\rm P.}\ref{prop:meet-neg}~(1a)}{=}
x\sqcap \neg \neg y 
\overset{{\rm P.}\ref{prop:axiom1a}~(1a)}{=} x\sqcap (y\sqcap y) \overset{{\rm T.}\ref{thm:firstindeax}~(1a)}{=} x\sqcap y.
\end{aligned}$

\noindent (3a)\quad 
$\begin{aligned}
x\sqcap\neg (\neg x\sqcap y) & \overset{{\rm T.}\ref{thm:firstindeax}~(1a)}{=}  x\sqcap\neg (\neg x\sqcap (y\sqcap y)) 
\overset{{\rm P.}\ref{prop:axiom1a}~(6a)}{=}  x\sqcap\neg (\neg x\sqcap \neg\neg y) \\
&\overset{{\rm def.}\vee}{=} x \sqcap (x \vee \neg y) \overset{{\rm T.}\ref{thm:axiom6ad}~(1a)}{=}  x\sqcap x. 
\end{aligned}$



\noindent (4a)\quad 
$\begin{aligned}
\neg x\sqcap\neg (x\sqcap y)
&\overset{{\rm P.}\ref{prop: De Morgan}~(1a)}{=} \neg x \sqcap (\neg x \vee \neg y)
\overset{{\rm T.}\ref{thm:axiom6ad}~(1a)}{=} \neg x \sqcap \neg x \overset{{\rm P.}\ref{prop:axiom1a}~(6a)}{=} \neg x
\end{aligned}$

\noindent (5a) \quad
$\begin{aligned}[t]
&\phantom{\overset{{\rm D.\ref{def:D-core+} (000)}}{=}} \neg x\sqcap \neg((x\sqcap y)\sqcap z)
\overset{{\rm P.}\ref{prop: De Morgan}(1a)}{=}
\neg x\sqcap (\neg(x\sqcap y) \vee \neg z)\\
&\overset{{\rm D.}\ref{def:D-core+}(4a)}{=}
(\neg x\sqcap \neg(x\sqcap y)) \vee (\neg x \sqcap \neg z)
\overset{{\rm P.}\ref{prop:meet-neg}(4a)}{=}
\neg x\vee (\neg x\sqcap \neg z) \\
&\overset{{\rm  P.}\ref{prop: De Morgan}(2a)}{=}
\neg x\vee \neg (x\vee z) 
\overset{{\rm  P.}\ref{prop: De Morgan}(1a)}{=}
\neg (x\sqcap (x\vee z)) 
\overset{{\rm T.}\ref{thm:axiom6ad}(1a)}{=}
\neg (x\sqcap x) \\
&\overset{{\rm D.}\ref{def:D-core+}(2a)}{=} \neg x
\end{aligned}$

\noindent (6a) \quad
$\begin{aligned}[t]
x\sqcap y
&\overset{\text{P.\,\ref{prop:axiom1a}(2a)}}{=}
(x\sqcap y)\sqcap (x\sqcap y)
\overset{\text{T.\,\ref{thm:axiom6ad}(1a)}}{=}
(x\sqcap y)\sqcap ((x\sqcap y)\vee (x\sqcap z))\\
&\overset{\rm{D.} \ref{def:D-core+}(4a)}{=}(x\sqcap y)\sqcap (x\sqcap (y\vee z))
\end{aligned}$.  

Setting $x:=\neg (y\vee z)$ we get

$\begin{aligned}[t]
&\neg(y\vee z)\sqcap y = (\neg(y\vee z)\sqcap y)\sqcap (\neg(y\vee z)\sqcap (y\vee z))
\overset{\rm{D.} \ref{def:D-core+}(4a)}{=} 
(\neg(y\vee z)\sqcap y)\sqcap \bot \\
&\overset{{\rm P.} \ref{proaxiom6}(4a)}{=} \bot. \qquad \text{i.e.}\quad \bot = \neg(y\vee z)\sqcap y \overset{{\rm P.}\ref{prop: De Morgan}~(2a)}{=} (\neg y\sqcap \neg z)\sqcap y  . 
\end{aligned}$. 

Setting $y:=x$ and $z:=\neg y$ we get 

$
\begin{aligned}
    \bot = (\neg x \sqcap \neg\neg y)\sqcap x \overset{{\rm P.}\ref{prop:axiom1a}~(1a)}{=}
    (\neg x \sqcap (y\sqcap y))\sqcap x 
    \overset{{\rm T.}\ref{thm:firstindeax}~(1a)}{=} (\neg x \sqcap y)\sqcap x 
\end{aligned}
$

        
\end{proof}
\begin{proposition}
\label{pro4}
     {\rm Let $\textbf{D}$ be a D-core algebra and  $x, y,z \in D$. We have:
 
    \noindent  $\begin{array}{ll}
    (1a)~((x\sqcap y)\sqcap z)\sqcap \neg x=\bot& (1b)~((x\sqcup y)\sqcup z)\sqcup \lrcorner x=\top\\
    (2a)~\neg(x\sqcap y)\sqcap \neg(x\sqcap\neg y)=\neg x&(2b)~\lrcorner(x\sqcup y)\sqcup \lrcorner(\lrcorner y\sqcup x)=\lrcorner x\\
    (3a)~\neg(\neg (x\sqcap (y\sqcap z))\sqcap z)\sqcap  z=x\sqcap(y\sqcap z)&\\
    (3b)~\lrcorner(\lrcorner (x\sqcup (y\sqcup z))\sqcup z)\sqcup z= x\sqcup(y\sqcup z)&
    \end{array}$}
\end{proposition}

\begin{proof}
    Let $x, y,z \in D$





        
\noindent (1a) \quad $\begin{aligned}[t]
\bot
&\overset{{\rm P.}\ref{prop:meet-neg}(6a)}{=}
((x\sqcap y)\sqcap z)\sqcap (\neg x \sqcap \neg ((x\sqcap y)\sqcap z))
\overset{\text{P.\,\ref{prop:meet-neg}(5a)}}{=}
((x\sqcap y)\sqcap z)\sqcap \neg x
\end{aligned}$

   \noindent (2a)\quad
$\begin{aligned}[t]
&\phantom{\overset{{\rm D.\ref{def:D-core+} (000)}}{=}} \neg(x\sqcap y)\sqcap \neg(x \sqcap \neg y)\overset{{\rm P.}\ref{prop: De Morgan}~(2a)}{=}
\neg((x\sqcap y)\vee (x\sqcap \neg y)) \\
&\overset{{\rm D.}\ref{def:D-core+}~(4a)}{=} \neg (x\sqcap (y\vee\neg y)) 
\overset{{\rm def.}\vee}{=}
\neg (x\sqcap \neg(\neg y\sqcap\neg\neg y))\overset{\text{D.\,\ref{def:D-core+}\,(6a)}}{=}
\neg(x\sqcap \neg\bot)\\
&\overset{{\rm T.}\ref{thm:firstindeax}~(3a)}{=}
\neg(x\sqcap (\top\sqcap\top))
\overset{\text{T.\,\ref{thm:firstindeax}\,(1a)}}{=}
\neg(x\sqcap \top)
\overset{{\rm P.}\ref{prop:axion6}~(1a)}{=}
\neg(x\sqcap x) \overset{{\rm D.}\ref{def:D-core+}~(2a)}{=}
\neg x
\end{aligned}$

        

   	
\noindent (3a)\quad $\begin{aligned}[t]
x\sqcap (y\sqcap z)
&\overset{{\rm D.}\ref{def:D-core+}~(5a)}{=}\neg\neg (x\sqcap (y\sqcap z))\\
&\overset{{\rm P.}\ref{pro4}(2a)}{=}
\neg(\neg(x\sqcap (y\sqcap z))\sqcap z)\sqcap \neg(\neg (x\sqcap (y\sqcap z))\sqcap \neg z)\\
&\overset{{\rm P.}\ref{prop:meet-neg}(5a)}{=}
\neg (\neg(x\sqcap(y\sqcap z))\sqcap z)\sqcap \neg\neg z\\
&\overset{{\rm P.}\ref{prop:axiom1a}~(1a)}{=}
\neg (\neg(x\sqcap(y\sqcap z))\sqcap z)\sqcap (z\sqcap z)\\
&\overset{{\rm T.}\ref{thm:firstindeax}~(1a)}{=}
\neg (\neg(x\sqcap(y\sqcap z))\sqcap z)\sqcap z
\end{aligned}$

    	

   
\end{proof}

\begin{proposition}
   \label{axiomassio} {\rm Let $\textbf{D}$ be a D-core algebra and  $x, y,z \in D$. Then, the following hold.\\
    \noindent  $\begin{array}{ll}       
    (1a)~\neg (x\sqcap \neg(y\sqcap z))=\neg(x\sqcap\neg y)\sqcap\neg(x\sqcap\neg z) &\\(1b)~\lrcorner(x\sqcup \lrcorner(y\sqcup z))=\lrcorner(x\sqcup\lrcorner y)\sqcup\lrcorner(x\sqcup\lrcorner z)&\\
     (2a)~x\sqcap (\neg(x\sqcap \neg y)\sqcap\neg(x\sqcap \neg z))=(y\sqcap x)\sqcap z & \\(2b)~x\sqcup (\lrcorner(x\sqcup \lrcorner y)\sqcup\lrcorner(x\sqcup \lrcorner z))=(y\sqcup x)\sqcup z.
    \end{array}$}
\end{proposition}
\begin{proof}
  Let $x, y,z \in D$

\noindent (1a)\quad 
$\begin{aligned}[t]
\neg(x\sqcap \neg(y\sqcap z)) &\overset{\text{P\ref{prop: De Morgan}~(1a)}}{=} \neg(x\sqcap (\neg y\vee \neg z)) \overset{\text{D\ref{def:D-core+}~(4a)}}{=} \neg ((x\sqcap \neg y) \vee (x\sqcap \neg z)) \\
&\overset{\text{P\ref{prop: De Morgan}~(2a)}}{=}  \neg (x\sqcap \neg y) \sqcap \neg (x\sqcap \neg z). 
\end{aligned}$

\noindent (2a)\quad 
$\begin{aligned}[t]
&\phantom{\overset{{\rm D.\ref{def:D-core+} (000)}}{=}} (y\sqcap x)\sqcap z
\overset{\text{D\ref{def:D-core+}~(1a)}}{=}
   z\sqcap (y\sqcap x)
   \overset{\text{P\ref{pro4} (3a)}}{=}
   \neg(\neg(z\sqcap (y\sqcap x))\sqcap x)\sqcap x
\\
&\overset{\text{D\ref{def:D-core+} (1a)}}{=}
   \neg(x\sqcap \neg((y\sqcap x)\sqcap z))\sqcap  x
\overset{\text{P\ref{axiomassio} (1a)}}{=}
   (\neg(x\sqcap\neg(y\sqcap x))\sqcap \neg(x\sqcap\neg z))\sqcap  x
\\
&\overset{{\rm P.}\ref{prop:meet-neg}~(1a)}{=}
   (\neg(x\sqcap\neg y)\sqcap \neg(x\sqcap\neg z))\sqcap x
\overset{\text{D\ref{def:D-core+} (1a)}}{=}
   x\sqcap (\neg(x\sqcap\neg y)\sqcap \neg(x\sqcap\neg z))
\end{aligned}
$
        \end{proof}
The following theorem shows that the axioms $(10a)$ and $(10b)$ in Definition~\ref{def:DBA} of a dBa can be derived from the axioms of a D-core algebra.
\begin{theorem}
\label{derivable axioms}
    {\rm Let $\textbf{D}$ be a D-core algebra and  $x, y,z \in D$. Then, the following hold.\\
    $\begin{array}{ll}
    (a)~x \sqcap ( y \sqcap  z) = (x \sqcap  y) \sqcap  z& (b)~x \sqcup (y \sqcup  z) = (x \sqcup  y)\sqcup  z.
    \end{array}$}
\end{theorem}
\begin{proof} 
Let $x,y, z\in D$

 \noindent
$\begin{aligned}[t]
& \phantom{\overset{{\rm D.\ref{def:D-core+} (000)}}{=}} x\sqcap (y\sqcap z)
\overset{\text{D\ref{def:D-core+} (1a)}}{=}
   x\sqcap (z\sqcap y) 
\overset{\text{P\ref{pro4} (3a)}}{=}
   \neg(\neg(x\sqcap (z\sqcap y))\sqcap y)\sqcap y
\\
& \overset{\text{D\ref{def:D-core+} (1a)}}{=}
   \neg(y\sqcap \neg(x\sqcap (z\sqcap y)))\sqcap y 
\overset{\text{P\ref{axiomassio} (1a)}}{=}
    (\neg(y\sqcap\neg x)\sqcap \neg(y\sqcap\neg(z\sqcap y)))\sqcap y
    \\&\overset{\text{P\ref{prop:meet-neg} (1a)}}{=}(\neg(y\sqcap\neg x)\sqcap\neg (y\sqcap\neg z))\sqcap y
\overset{\text{D\ref{def:D-core+} (1a)}}{=}
    y\sqcap (\neg(y\sqcap\neg x)\sqcap \neg(y\sqcap\neg z))\\
&\overset{\text{P\ref{axiomassio} (2a)}}{=}
   (x\sqcap y)\sqcap z
\end{aligned}$

\end{proof}
Summarizing the results so far, we obtain the following main theorem. It will be used in the next section to refine representation theorems for dBas.

\begin{theorem}
    \label{D+coreanddba}{\rm An 
    algebra $\textbf{D}:=(D, \sqcap,\sqcup, \neg, \lrcorner, \top, \bot)$ of type (2,2,1,1,0,0) is a dBa if and only if it is a D-core algebra.}
\end{theorem}
\begin{proof}
    The proof follows from Theorems \ref{thm:firstindeax}, \ref{thm:axiom6ad} and 
    \ref{derivable axioms}.
\end{proof}

It can be shown that Definition~\ref{def:D-core+} represents a minimal set of axioms. 
To demonstrate this, we must verify the independence of each pair of axioms from the remaining ones. 
The full proof is quite tedious. Therefore, we illustrate it with the proof for the pair (5a),(5b) only.

\begin{theorem}
\label{cexample7}
    {\rm There is an algebra $\textbf{D}:=(D, \sqcap,\sqcup, \neg, \lrcorner, \top, \bot)$ that satisfies $(1a)-(4a)$, $(1b)-(4b)$, $(6a), (6b)$ and $(7)$ and does not satisfy $(5a)$ and $(5b)$.}
\end{theorem}
\begin{proof}
In this proof, we construct a two-element  algebra in which axioms (5a) and (5b) do not hold. Let $\mathbf{D} := (\{a, b\}; \sqcap, \sqcup, \neg, \lrcorner, \top, \bot)$, where $a=\bot$, $b=\top$ and the operations are defined by the following tables:

\begin{center}
\begin{tabular}{|c|c|c|}
\hline
$\sqcap$ & a & b \\
\hline
a & a & a \\
b & a & b \\
\hline
\end{tabular}
\quad
\begin{tabular}{|c|c|c|}
\hline
$\sqcup$ & a & b \\
\hline
a & a & b \\
b & b & b \\
\hline
\end{tabular}
\quad
\begin{tabular}{|c|c|c|}
\hline
$\wedge$ & a & b \\
\hline
a & b & b \\
b & b & b \\
\hline
\end{tabular}
\quad
%
\begin{tabular}{|c|c|c|}
\hline
$\vee$ & a & b \\
\hline
a & a & a \\
b & a & a \\
\hline
\end{tabular}
\quad
\begin{tabular}{|c|c|c|}
\hline
$x$ & $\neg x$ & $\lrcorner x$\\
\hline
a & a & b\\
b & a & b\\
\hline
\end{tabular}
\end{center}
Now from the above tables, we get 
$\neg\neg(b\sqcap b) = \neg\neg b = \neg a = a \neq b = b\sqcap b$, failing (5a). Similarly $\lrcorner\lrcorner(a\sqcup a) = b \neq a = a\sqcup a$, failing (5b).
%
%
\end{proof}
\section{Refinement of the Representation Theorem for dBas}
\label{sec:re} In this section, we build on the refined framework of double Boolean algebras to revisit existing results and derive more compact formulations, leading to stronger and clearer theorems.

 \subsection{Boolean Representations of dBas}

We still use $\mathbf{ P}$  and $\mathbf{Q}$ to denote the two Boolean algebras  defined in Section~\ref{pre}.  Recall that Theorem \ref{dbaboolean} shows that a dBa {\bf A} can be constructed from them with two pairs of maps $(r,e)$ and $(r',e')$ satisfying three conditions: the first ensures Axiom~12, the second ensures Axioms~(4a) and (4b), and the third ensures Axiom~(11a)  and (11b). Now, because of Theorems \ref{thm:firstindeax} and \ref{D+coreanddba}, we can see that the third condition is actually redundant. Hence, we obtain the following result.

\begin{theorem}
 	\label{newdbaboolean}{\rm $\textbf{A}$ is a dBa iff the following holds
 		\begin{enumerate}
 		\item $e\circ r\circ e^{\prime}\circ r^{\prime}=e^{\prime}\circ r^{\prime}\circ e\circ r$
 		\item $e(r(x)\wedge_p r(e^{\prime}(r^{\prime}(x)\vee_q r^{\prime}(y))))=e(r(x))$ and $e^{\prime}(r^{\prime}(x)\vee_q r^{\prime}(e(r(x)\wedge_p r(y))))=e^{\prime}(r^{\prime}(x))$ for all $x, y\in A$.
 		\end{enumerate}
 		Moreover every dBa can be obtained from such a construction.}
 \end{theorem}
In the next theorem, we establish a more general result extending Theorem~\ref{thm:pure-trivial-dBa}. For that, we define the \emph{generalized linear sum} of two posets $(P, \leq_{P})$ and $(Q, \leq_{Q})$, denoted by $P \oplus_{g} Q$, in the same way as the ordinary linear sum, except that  $P$ and $Q$ may be not disjoint.  

\begin{definition}\label{gls}
Formally, the underlying set of $P \oplus_{g} Q$ is $P \cup Q$, and the 
relation $\leq_{P \oplus_{g} Q}$ on $P \oplus_{g} Q$ is given by
\[
x \leq_{P \oplus_{g} Q} y :\iff
\begin{cases}
x, y \in P \text{ and } x \leq_{P} y, & \text{or} \\[4pt]
x, y \in Q \text{ and } x \leq_{Q} y, & \text{or} \\[4pt]
x \in P \text{ and } y \in Q & \text{or} \\[4pt]
x=\bot_{q} ~\text{and}~ y=\top_{p}
\end{cases}
\]
\end{definition}

That is, every element of $P$ is declared less than every element of $Q$, but  $P \cap Q$ may be nonempty.

\begin{proposition}
    {\rm Let $(P, \leq_{P})$ and $(Q, \leq_{Q})$ be two posets then $P \oplus_{g} Q$ is a quasiordered set.}
\end{proposition}
\begin{proof}
By definition, the relation $\le_{P \oplus_g Q}$ is reflexive. 
Let $x \le_{P \oplus_g Q} y$ and $y \le_{P \oplus_g Q} z$. 
If $x, y, z \in P$ or $x, y, z \in Q$, then transitivity follows from that of $\le_P$ and $\le_Q$, respectively. 
If one of $y$ or $z$ lies in $Q$ while $x \in P$, then by the definition of $\le_{P \oplus_g Q}$ we also have $x \le_{P \oplus_g Q} z$. 
Hence, $\le_{P \oplus_g Q}$ is transitive.
\end{proof}

The relation $\le_{P \oplus_g Q}$ is not necessarily antisymmetric. 
For example, choose 
$x, y \in P \cap Q$ with $x \neq y$. 
Then, both $x \le_{P \oplus_g Q} y$ and $y \le_{P \oplus_g Q} x$ hold. 
Moreover, when $P \cap Q = \emptyset$, 
then the generalized linear sum reduces to the ordinary linear sum.

For the Boolean algebras {\bf P} and {\bf Q}, let us consider the pairs of maps $r:P\cup Q\rightleftharpoons P: e$ and $r^{\prime}:P\cup Q \rightleftharpoons Q: e^{\prime}$ defined by
\begin{center}
 \begin{tabular}{ll}
$r(x)=\begin{cases}
			x, & \text{if $x\in P$}\\
           \top_{p}, & \text{otherwise}
		 \end{cases}$   &  $r^{\prime}(x)=\begin{cases}
			x, & \text{if $x\in Q$}\\
           \bot_{q}, & \text{otherwise}
		 \end{cases}$\\
$e(x)=x$ for all $x\in P$    & $e^{\prime}(x)=x$ for all $x\in Q$
\end{tabular}   
\end{center}
Then, we can define the following universal algebra,
        \[ \textbf{P+Q}:=(P\cup Q, \sqcap, \sqcup, \neg,\lrcorner, e^{\prime}(\top_{q}), e(\bot_{p}))\] 
where the operatons are defined as in Equations~(\ref{eqn:dba_ops_from_e-r-pairs}).  
    In addition, we define a relation $\sqsubseteq$ on $P\cup Q$, as in Equation~(\ref{Eq:dBa quasiorder}).
    

      \begin{lemma}
         \label{dBaquasi} {\rm For all $x, y\in P\cup Q$, $x\sqsubseteq y$ if and only if $x\leq_{P \oplus_{g} Q} y$.}
      \end{lemma}
      \begin{proof}
          Let $x, y\in P\cup Q$ and $x\sqsubseteq y$ then $x\sqcap y=e(r(x)\wedge_{p} r(x))$ and $x\sqcup y=e^{\prime}(r^{\prime}(y)\vee_{q} r^{\prime}(y))$. If $x, y\in P$ then $x\wedge_{p} y=e(r(x)\wedge_{p} r(y))=x\sqcap y=e(r(x)\wedge_{p} r(x))=x\wedge_{p}x=x$ which implies that $x\leq_{P} y$. So $x\leq_{P \oplus_{g} Q} y$.  If $x, y\in Q$, then $x\vee_{q} y=e^{\prime}(r^{\prime}(x)\vee_{q} r^{\prime}(y))=x\sqcup y=e^{\prime}(r^{\prime}(y)\vee_{q} r^{\prime}(y))=y\vee_{q}y=y$, which implies that $x\leq_{Q} y$. So $x\leq_{P \oplus_{g} Q} y$. If $x\in P$ and $y\in Q$, then $x\leq_{P \oplus_{g} Q} y$.  Now, if $x\notin P$ and $y\in P$, $e(r(x)\wedge_{p} r(y))=e(r(x)\wedge_{p} r(x))$ implies that $y=\top_{p}$ and $e^{\prime}(r^{\prime}(x)\vee_{p} r^{\prime}(y))=e^{\prime}(r^{\prime}(y)\vee_{p} r^{\prime}(y))$, implies that $x=\bot_{q}$. So $x\leq_{P \oplus_{g} Q} y$.

          Conversely, let $x, y\in P\cup Q$ and $x\leq_{P \oplus_{g} Q} y$. If $x, y\in P$ then $x\sqcap y=x\wedge_{p} y=x\wedge_{p} x=x\sqcap x$ and $x\sqcup y=e^{\prime}(r^{\prime}(x)\vee_{q} r^{\prime}(y))=e^{\prime}(\bot_{q}\vee\bot_{q})=y\sqcup y$. So $x\sqsubseteq y$. If $x, y\in Q$ then we can show that $x\sqcap y=e(\top_{p}\wedge_{p}\top_{p})=x\sqcap x$ and $x\sqcup y=x\vee_{q} y=y\vee_{q} y=y\sqcup y$, which implies that $x\sqsubseteq y$.  If $x\in P$ and $y\in Q$ then $x\sqcap y=x\wedge_{p} x=x\sqcap x$ and $x\sqcup y=y\vee_{q} y=y\sqcap y$ which implies that $x\sqsubseteq y$. If $x=\bot_{q}$ and $y=\top_{p}$ then $x\sqcap y=\top_{p}=x\sqcap x$ and $x\sqcup y=\bot_{q}=y\sqcup y$ which implies that $x\sqsubseteq y$.
          
      \end{proof}

      \begin{definition}\label{def:generalized glued sum}
          {\rm A generalized glued sum of {\bf P} and {\bf Q} is  an ordered algebraic structure $\textbf{P$\oplus$ Q}:=(P\cup Q, \sqcap, \sqcup, \neg,\lrcorner, \top_{q}, \bot_{p})$ such that.
          \begin{itemize}
              \item[(a)] For $x, y\in P$, $x\sqcap y=x\wedge_{p} y$, $\neg x=\neg_{p}x$
              \item[(b)] For $x, y\in Q$, $x\sqcup y=x\vee_{q} y$, $\lrcorner x=\neg_{q} x$
              \item[(c)]  $x\leq_{P \oplus_{g} Q} y$ if and only if $x\sqcap y=x\sqcap x~\mbox{and}~x\sqcup y=y\sqcup y~\mbox{for all}~x, y\in P\cup Q$
          \end{itemize}}
      \end{definition}

      \begin{proposition}\label{prop:PQ-generalized glued sum}
          {\rm  $\textbf{P+Q}:=(P\cup Q, \sqcap, \sqcup, \neg,\lrcorner, e^{\prime}(\top_{q}), e(\bot_{p}))$ is a  generalized glued sum.}
      \end{proposition}
       \begin{proof}
           Condition (a) and (b) follows from the definition of the operations. Condition (c) follows from Lemma \ref{dBaquasi}. 
       \end{proof}

       Before going to the next result we define generalized D-core algebra.
       \begin{definition}
           {\rm An universal algebra $(D, \sqcap,\sqcup, \neg,\lrcorner, \top, \bot)$ is said to be {\it generalized D-core algebra} if it satisfies (1a)-(2a), (4a)-(6a), (1b)-(2b), (4b)-(6b) and 7 of Definition \ref{def:D-core+}.}
       \end{definition}

\begin{example}

  { \rm  In this example we produced an universal algebra that is a generalized D-core algebra but not D-core.    
  
  \noindent $\begin{array}{|c|ccc|}
  \hline
\sqcup & a & b & c \\ \hline
a & c & b & c \\
b & b & b & b \\
c & c & b & c\\
\hline
\end{array}
\quad
\begin{array}{|c|ccc|}
\hline
\wedge & a & b & c \\ \hline
a & c & c & c \\
b & c & b & c \\
c & c & c & c\\\hline
\end{array}
\quad
\begin{array}{|c|ccc|}
\hline
\sqcap & a & b & c \\ \hline
a & a & c & b \\
b & c & b & a \\
c & b & a & c\\\hline
\end{array}
\quad
 \begin{array}{|c|ccc|}
 \hline
\vee & a & b & c \\ \hline
a & a & c & b \\
b & c & b & a \\
c & b & a & c\\\hline
\end{array}
\quad
 \begin{array}{|c|cc|}
 \hline
x & \neg x & \lrcorner x \\ \hline
a & a & b  \\
b & c & c  \\
c & b & b \\\hline
\end{array}
$

\noindent $a\sqcap (a\sqcup c)=b\neq a=a\sqcap a$ and $a\sqcup (a\sqcap c)=b\neq c=a\sqcup a$}
\end{example}
The above example implies that there is a non D-core algebra $\textbf{D}$ that satisfies $(1a)-(2a)$, $(1b)-(2b)$, $(4a)-(6a), (4b)-(6b)$ and $(7)$ and does not satisfies $(3a)$ and $(3b)$.
\begin{theorem}
\label{consgeneralizedD}	{\rm For the two Boolean algebras {\bf P} and {\bf Q}, assume that $\{\top_p,\bot_q\}\subseteq P\cap Q$.  Then, $\textbf{P}+\textbf{Q}$ is a  generalized D-core algebra.}
\end{theorem}
\begin{proof}
We have to show that $e\circ r\circ e^{\prime}\circ r^{\prime}(x)=e^{\prime}\circ r^{\prime}\circ e\circ r(x)$ for all $x\in P\cup Q$. For that, we have to consider three cases:
\begin{description}
\item $x\in P\cap Q$:  then  $e\circ r\circ e^{\prime}\circ r^{\prime}(x)=x=e^{\prime}\circ r^{\prime}\circ e\circ r(x)$ 
\item  $x\in P$ and $x\notin Q$: then $e\circ r\circ e^{\prime}\circ r^{\prime}(x)=e\circ r(\bot_q)=\bot_q$ as $\bot_q\in P$ and $e^{\prime}\circ r^{\prime}\circ e\circ r(x)=e^{\prime}\circ r^{\prime}(x)=\bot_{q}$ as $x\notin Q$. 
\item $x\notin P$ and $x\in Q$: then $e\circ r\circ e^{\prime}\circ r^{\prime}(x)=e\circ r(x)=\top_p$ as $x\notin P$  and $e^{\prime}\circ r^{\prime}\circ e\circ r(x)=e^{\prime}\circ r^{\prime}(\top_p)=\top_p$ as $\top_p\in Q$.
\end{description}
\end{proof}


   	
\begin{corollary}
{\rm Assume that the domains of the  Boolean algebras {\bf P} and {\bf Q} satisfy $P\cap Q=\{\top_{p}\}=\{\bot_{q}\}$, then $\mathbf{P} + \mathbf{Q}$ is a pure and trivial double Boolean algebra. Conversely, every pure and trivial double Boolean algebra is glued sum of two Boolean algebras.}
\end{corollary}

\begin{proof}

By Theorem~\ref{consgeneralizedD}, the algebra $\mathbf{P}+\mathbf{Q}$ is a generalized $D$-core algebra. For any $x,y\in P\cup Q$, we have
\[
e\!\left(r(x)\wedge_p r\!\left(e'(r'(x)\vee_q r'(y))\right)\right)
= x\sqcap (x\sqcup y)
\quad\text{and}\quad
e(r(x))=x\sqcap x.
\]
Thus, in order to show that $\mathbf{P}+\mathbf{Q}$ is a $D$-core algebra, it suffices to verify that $x\sqcap (x\sqcup y)=x\sqcap x$
for all  $x,y\in P\cup Q$.

This equality holds in all possible cases determined by the membership of $x$ and $y$. If $x\notin P$ and $y\notin Q$, then $x\sqcap (x\sqcup y)=x\sqcap x$. If $x\notin Q$ and $y\notin P$, then $x\sqcap (x\sqcup y)=x\sqcap y=x=x\sqcap x$. If $x,y\in P$ and $x,y\notin Q$, then $x\sqcup y=\bot_q$, and hence $x\sqcap (x\sqcup y)=x\sqcap\bot_q=x=x\sqcap x$. If $x,y\in Q$ and $x,y\notin P$, then $x\sqcup y=x\vee_q y$, and either $x\vee_q y\notin P$, in which case $x\sqcap (x\sqcup y)=\top_p=x\sqcap x$, or $x\vee_q y\in P\cap Q$, which implies $x\vee_q y=\top_p=\bot_q$ and again $x\sqcap (x\sqcup y)=x\sqcap\top_p=\top_p=x\sqcap x$. Finally, when one or both of the elements equals the common boundary element $\top_p=\bot_q$, a direct verification shows that $x\sqcap (x\sqcup y)=x\sqcap x$ in every remaining sub-case. Consequently, $x\sqcap (x\sqcup y)=x\sqcap x$ for all $x,y\in P\cup Q$, and hence $\mathbf{P}+\mathbf{Q}$ is a $D$-core algebra. By Theorem~\ref{newdbaboolean}, it follows that $\mathbf{P}+\mathbf{Q}$ is a dBa. Moreover, for every $x\in P\cup Q$, we have $x\sqcap x=x$ when $x\in P$ and $x\sqcup x=x$ when $x\in Q$, showing that $\mathbf{P}+\mathbf{Q}$ is a pure dBa. Finally, since $\top_q\sqcap\top_q=\top_p=\bot_q=\bot_p\sqcup\bot_p$, the algebra $\mathbf{P}+\mathbf{Q}$ is a trivial dBa.

       Conversely, let $\textbf{D}:= (D;  \sqcap, \sqcup,\neg,\lrcorner,\top,\bot)$ be a pure and trivial dBa. For the Boolean algebras  $\textbf{P}=(D_{\sqcap}, \sqcap,\vee,\neg, \neg\bot, \bot)$ and $\textbf{Q}=(D_{\sqcup}, \sqcup, \wedge, \lrcorner, \top, \lrcorner\top)$, let us  define the maps $r: D_{\sqcup}\cup D_{\sqcap}\rightarrow D_{\sqcap}$ by $r(x)=x\sqcap x$, $e:D_{\sqcap}\rightarrow  D_{\sqcup}\cup D_{\sqcap}$ by $e(x)=x$, $r^{\prime}:D_{\sqcup}\cup D_{\sqcap}\rightarrow D_{\sqcup}$  by $r^{\prime}(x)=x\sqcup x$, and  $ e^{\prime}: D_{\sqcup}\rightarrow D_{\sqcup}\cup D_{\sqcap}$  by $e^{\prime}(x)=x$.  Then, by Theorem \ref{newdbaboolean}, $(D_{\sqcap}\cup D_{\sqcup}, \sqcap_{r}, \sqcup_{r^{\prime}}, \neg_{r},\lrcorner_{r^{\prime}}, e^{\prime}(\top), e(\bot))$ is a  dBa, where, $x\sqcap_{r} y:=e(r(x)\sqcap r(y)),~ x\sqcup_{r^{\prime}} y:=e^{\prime}(r^{\prime}(x)\sqcup r^{\prime}(y)),~ \neg_{r} x:=e(\neg r(x))$, and $\lrcorner_{r} x:=e^{\prime}(\lrcorner r^{\prime}(x))$.  Moreover, for all $x, y\in D_{\sqcap}\cup D_{\sqcup}$, $x\sqcap_{r} y=x\sqcap y,~ x\sqcup_{r^{\prime}} y=x\sqcup y,~ \neg_{r} x=\neg x$, and $\lrcorner_{r} x=\lrcorner x$, $e^{\prime}(\top)=\top$ and $e(\bot)=\bot$. To complete the proof we need to show that  $x\leq_{D_{\sqcap} \oplus_{g} D_{\sqcup}} y$ if and only if $x\sqcap y=x\sqcap x~\mbox{and}~x\sqcup y=y\sqcup y~\mbox{for all}~x, y\in D_{\sqcap}\cup D_{\sqcup}$.

        For $x, y\in D_{\sqcap}$, $x\leq_{D_{\sqcap} \oplus_{g} D_{\sqcup}} y$ if and only if $x\sqsubseteq_{D_\sqcap} y$ if and only if $x\sqsubseteq y $ if and only if $x\sqcap y=x\sqcap x~\mbox{and}~x\sqcup y=y\sqcup y$

       Similar to the above, we can show that $x, y\in D_{\sqcup}$, $x\leq_{D_{\sqcap} \oplus_{g} D_{\sqcup}} y$ if and only if $x\sqcap y=x\sqcap x~\mbox{and}~x\sqcup y=y\sqcup y$.

        Now if $x\in D_{\sqcap}$,  $y\in D_{\sqcup}$ and  $x\sqcap y=x\sqcap x~\mbox{and}~x\sqcup y=y\sqcup y$ then $x\leq_{D_{\sqcap} \oplus_{g} D_{\sqcup}} y$ by definition. Now we assume that $x\leq_{D_{\sqcap} \oplus_{g} D_{\sqcup}} y$ and  $x\in D_{\sqcap}$, $y\in D_{\sqcup}$. $x\sqcap y=x\sqcap (y\sqcap y)= x\sqcap (\top\sqcap \top)=x\sqcap \top=x\sqcap x$ and $x\sqcup y=(x\sqcup x) \sqcap y
        =(\bot\sqcup \bot)\sqcup y
        =\bot\sqcup y
        =y\sqcup y$.
         
        

\end{proof}





\subsection{Representation Theorems of dBAs}\label{see:tre}
\noindent \\
In this subsection we will discussed topological representation of  dBas.
Recall the topological spaces  $(\mathcal{F}_{pr}(\textbf{D}),\mathcal{T})$ and $(\mathcal{I}_{pr}(\textbf{D}),\mathcal{J})$ defined in Section \ref{pre}. By Proposition \ref{protopo} these topological spaces are  compact and totally disconnected topological spaces. Moreover, we can prove the following. 

\begin{lemma} 
	\label{clopenset}	{\rm A subset $X$ (resp. $Y$) is clopen in $\mathcal{F}_{pr}(\textbf{D})$ (resp. $\mathcal{I}_{pr}(\textbf{D})$) if and only if $X=F_{x}$ (resp. $Y=I_{x}$) for some $x\in D$.}
\end{lemma}
\begin{proof}
	Since $X$ is closed, we have $X = \bigcap_{j \in J} \bigcup_{a \in D_j} F_a$, where for each $j \in J$ the set $D_j$ is a finite subset of $D$. As $X$ is open, its complement $X^c = \bigcup_{j \in J} \bigcap_{a \in D_j} F_a^c$ is a closed subset of $(F_{pr}(D), \mathcal{T})$. Since $(F_{pr}(D), \mathcal{T})$ is compact, $X^c$ is compact, and hence there exists a finite subset $E \subseteq J$ such that $X^c = \bigcup_{j \in E} \bigcap_{a \in D_j} F_a^c$. Using Lemma \ref{derivation}(6), (3) and Proposition \ref{prop: De Morgan}(2a), we obtain $X^c = \bigcup_{j \in E} \bigcap_{a \in D_j} F_{\neg a} = \bigcup_{j \in E} F_{\sqcap_{a \in D_j} \neg a} = \bigcup_{j \in E} F_{\neg(\bigvee_{a \in D_j} a)}$. Taking complements yields $X = \bigcap_{j \in E} F_{\neg\neg(\bigvee_{a \in D_j} a)} = \bigcap_{j \in E} F_{\bigvee_{a \in D_j} a}$, since $\bigvee_{a \in D_j} a \in D_{\sqcap}$ for each $j \in E$. Therefore, $X = F_{\sqcap_{j \in E} (\bigvee_{a \in D_j} a)}$, and setting $x = \sqcap_{j \in E} (\bigvee_{a \in D_j} a)$ we conclude that $X = F_x$.

\end{proof}

Let $clopen(\mathcal{J})$ denote the set of all clopen subsets of $\mathcal{J}$ and 
$clopen(\mathcal{T})$ the set of all clopen subsets of $\mathcal{T}$. 
Let $\mathcal{D}\subseteq clopen(\mathcal{T})\times clopen(\mathcal{J})$, such that $\mathcal{D}:=\{(F_{ x}, I_{x})~\mid~x\in D\}$. Let $\mathcal{D}_{\sqcap}:=\{(F_{ x}, I_{x\sqcap x})~\mid~ x\in D\}$ and $\mathcal{D}_{\sqcup}:=\{(F_{(x\sqcup x)}, I_{x})~\mid~ x\in D\}$. Note that $\mathcal{D}_{\sqcap}, \mathcal{D}_{\sqcup}\subseteq \mathcal{D}$ because $F_x=F_{x\sqcap x}$ and $I_x=I_{x\sqcup x}$.  We can then prove the following.

\begin{lemma}
{\rm Let $\mathbf{D}$ be a dBa. Then the following hold:

\begin{enumerate}
    \item $(\mathcal{D}_{\sqcap},\vee_{B}, \wedge_{B},  \neg_{B}, 0)$ is a Boolean algebra, where:
    \begin{itemize}
        \item Join: $(F_{ x}, I_{x\sqcap x}) \vee_{B} (F_{ y}, I_{y\sqcap y}) = (F_{ x \vee  y}, I_{  x \vee  y})$
        \item Meet: $(F_{ x}, I_{x\sqcap x}) \wedge_{B} (F_{ y}, I_{y\sqcap y}) = (F_{ x \sqcap  y}, I_{  x \sqcap  y})$
        \item Complement: $\neg_{B} (F_{ x}, I_{x\sqcap x}) = (F_{\neg x}, I_{\neg x})$
         \item Bottom element: $0 = (F_\bot, I_{\bot})$
    \end{itemize}

    \item $(\mathcal{D}_{\sqcup},\vee_{B},\wedge_{B}, \lrcorner_{B}, 1)$ is a Boolean algebra, where:
    \begin{itemize}
        \item Join: $(F_{(x \sqcup x)}, I_x) \vee_{B} (F_{(y \sqcup y)}, I_y) = (F_{(x \sqcup y)}, I_{x \sqcup y})$
        \item Meet: $(F_{(x \sqcup x)}, I_x) \wedge_{B} (F_{(y \sqcup y)}, I_y) = (F_{(x \wedge y)}, I_{x \wedge y})$
        \item Complement: $\lrcorner_{B}(F_{x \sqcup x}, I_x) = (F_{\lrcorner x}, I_{\lrcorner x})$
        \item Top element: $1 = (F_{\top}, I_\top)$
    \end{itemize}
\end{enumerate}}
\end{lemma}
\begin{proof}
    The proof is straightforward. The key observation is that the connective \(\wedge_{B}\) mimics the operation \(\sqcap\) of the dBa \textbf{D}, and the connective \(\vee_{B}\) mimics the operation \(\sqcup\) of \textbf{D}. We know that for a dBa \textbf{D}, the structures \((D_{\sqcap}, \sqcap,\vee, \neg, \bot)\) and \((D_{\sqcup}, \sqcup,\wedge, \lrcorner, \top)\) each form a Boolean algebra.

\end{proof}

\noindent Now we define the pairs of maps $(r_0,e_0)$ and $(r_1,e_1)$ as follows:
%
%
\begin{align*}
  r: \mathcal{D}\rightarrow \mathcal{D}_{\sqcap},~(F_{ x}, I_{x})\mapsto (F_{ x}, I_{x\sqcap x});~~~& e: \mathcal{D}_{\sqcap}\rightarrow \mathcal{D},~(F_{ x}, I_{x\sqcap x})\mapsto (F_{x\sqcap x}, I_{x\sqcap x}); \\ 
r^{\prime}: \mathcal{D}\rightarrow \mathcal{D}_{\sqcup},~(F_{ x}, I_{x})\mapsto (F_{x\sqcup x}, I_x);~~& e^{\prime}: \mathcal{D}_{\sqcup}\rightarrow \mathcal{D},~(F_{x\sqcup x}, I_x)\mapsto (F_{x\sqcup x}, I_{x\sqcup x}).  
\end{align*}



\begin{proposition}
	\label{repdba}	{\rm  The pairs  $(r, e)$ and $(r^{\prime}, e^{\prime})$ of maps satisfies the condition of Theorem \ref{newdbaboolean}.}
\end{proposition}
\begin{proof}
	From the definition of the map it is clear that  $r\circ e$ and $r^{\prime}\circ e^{\prime}$ are identity on $\mathcal{D}_{\sqcap}$ and $\mathcal{D}_{\sqcup}$, respectively. Moreover $r,r^{\prime}$ are surjective and $e,e^{\prime}$ injective. 
Let $(F_{x}, I_{x}), (F_{y}, I_{y})\in \mathcal{D}$. 	

\noindent 
$\begin{aligned}[t]
         e\circ r\circ e^{\prime}\circ r^{\prime}((F_{x}, I_{x})) &\overset{}{=}e\circ r\circ e^{\prime}((F_{x\sqcup x}, I_{x}))\overset{}{=}e\circ r((F_{x\sqcup x}, I_{x\sqcup x}))\\&\overset{}{=}e((F_{x\sqcup x}, I_{(x\sqcup x)\sqcap (x\sqcup x)})) \overset{}{=}(F_{(x\sqcup x)\sqcap (x\sqcup x)}, I_{(x\sqcup x)\sqcap (x\sqcup x)})
     \end{aligned}$
     
     \noindent Similar to the above we can prove that
     
     \noindent $\begin{aligned}[t]
         e^{\prime}\circ r^{\prime}\circ e\circ r((F_{x}, I_{x})) &\overset{}{=}(F_{(x\sqcap x)\sqcup (x\sqcap x)}, I_{(x\sqcap x)\sqcup (x\sqcap x)})\\& 
         =(F_{(x\sqcup x)\sqcap (x\sqcup x)}, I_{(x\sqcup x)\sqcap (x\sqcup x)})\overset{}{=}e\circ r\circ e^{\prime}\circ r^{\prime}((F_{x}, I_{x}))
     \end{aligned}$


     \noindent $\begin{aligned}[t]
          r^{\prime}((F_{x}, I_{x}))\vee r^{\prime}((F_{y}, I_{y})) &\overset{}{=}(F_{x\sqcup x}, I_{x})\vee (F_{y\sqcup y}, I_{y}) \overset{}{=}(F_{x\sqcup y},  I_{x\sqcup y})\\&\overset{\text{P.\,\ref{prop:axiom1a}(2b)}}{=}e(F_{(x\sqcup y)\sqcup (x\sqcup y)},  I_{x\sqcup y}) \overset{}{} ~\text{which implies that }
     \end{aligned}$
  
   \noindent $\begin{aligned}[t]
         & e^{\prime}(r^{\prime}((F_{x}, I_{x}))\vee r^{\prime}((F_{y}, I_{y}))) \overset{}{=}(F_{(x\sqcup y)\sqcup (x\sqcup y)},  I_{(x\sqcup y)\sqcup (x\sqcup y)})
         =(F_{x\sqcup y},  I_{x\sqcup y})
     \end{aligned}$
    
      \noindent So $r(e^{\prime}(r^{\prime}((F_{x}, I_{x}))\vee r^{\prime}((F_{y}, I_{y}))))=(F_{x\sqcup y}, I_{(x\sqcup y)\sqcap (x\sqcup y)})$.
      
      \noindent Hence, $r(F_{x}, I_{x})\wedge r(e^{\prime}(r^{\prime}((F_{x}, I_{x}))\vee r^{\prime}((F_{y}, I_{y}))))=(F_{x\sqcap (x\sqcup y)}, I_{x\sqcap (x\sqcup y)} )$ 
 which implies that 
 
 \noindent $\begin{aligned}[t]
     &e(r(F_{x}, I_{x})\wedge r(e^{\prime}(r^{\prime}((F_{x}, I_{x}))\vee r^{\prime}((F_{y}, I_{y})))))\overset{}{=}(F_{x\sqcap (x\sqcup y)}, I_{x\sqcap (x\sqcup y)} )\\&\overset{\text{D.\, \ref{def:D-core+}(3a)}}{=}(F_{x\sqcap x}, I_{x\sqcap x})\overset{}{=} (F_{x}, I_{x\sqcap x})\overset{}{=}e\circ r((F_{x}, I_{x}))
 \end{aligned}$
 
\noindent The third equality holds as $F_{x\sqcap x}=F_{x}$. 
%
%
Similarly, we can show that 

$e^{\prime}(r^{\prime}(F_{x}, I_{x})\wedge r^{\prime}(e(r((F_{x}, I_{x}))\vee r((F_{y}, I_{y})))))= e^{\prime}\circ r^{\prime}((F_{x}, I_{x}))$.
\end{proof}

\begin{corollary}\label{constalgeb}
\rm $(\mathcal{D},\, \sqcap,\, \sqcup,\, \neg,\, \lrcorner,\, e^{\prime}(1),\, e(0))$ is
 a dBa, where the operations are defined by Equations~(\ref{eqn:dba_ops_from_e-r-pairs}).
\end{corollary}

\begin{proof}
	{It follows from Theorem \ref{newdbaboolean} and Proposition \ref{repdba}.}
\end{proof}

\begin{theorem}
	\label{new:rep}{\rm Every dBa $\textbf{D}$ is quasi-isomorphic to the dBa of clopen subsets of the Stone topological spaces. Moreover if the dBa \textbf{D} is contextual then the quasi-isomorphism  becomes isomorphism.}
\end{theorem}
\begin{proof}
	We consider the Stone  topological spaces $(\mathcal{F}_{pr}(\textbf{D}),\mathcal{T})$ and $(\mathcal{I}_{pr}(\textbf{D}),\mathcal{J})$. Then the product $(\mathcal{F}_{pr}(\textbf{D})\times \mathcal{I}_{pr}(\textbf{D}),\mathcal{T}\times \mathcal{J} )$  is also a Stone topological space. By Lemma \ref{clopenset},  $F_{x}\times I_{x}$ are clopen subset in the product topology. The clopen set $F_{x}\times I_{x}$ can be identified with the pair $(F_{x}, I_{x})$.  Then from the Corollary \ref{constalgeb}, we get the required algebra, $\mathcal{D}:=(\mathcal{D}, \sqcap, \sqcup, \neg,\lrcorner, e^{\prime}(1), e(0))$. The map  $h:\textbf{D}\rightarrow\mathcal{D} $ is defined by $h(x)=(F_{x}, I_{x})$. This map is surjective. By Theorem \ref{firstrepsentation},  $h$ is a dBa homomorphism and $x\sqsubseteq y $ if and only if $h(x)\sqsubseteq h(y)$.
	
	Now if $\textbf{D}$ is contextual then the map $h$ is also injective as the quasi-order becomes partial order. Hence in this case $h$ is an isomorphism.
\end{proof} 

Recall Definition \ref{cTS} about context on topological spaces.

\begin{definition}
    {\rm Let \(\mathbb{K}^{T} := ((G, \rho), (M, \tau), R)\) be a CTS and let $A\subseteq G$ and $B\subseteq M$. Then, a pair of sets $(A, B)$ is a {\it clopen protoconcept} (resp.\ {\it semiconcept})  if $(A, B)$ is a protoconcept (resp.\ semiconcept) of $\mathbb{K}=(G, M, R)$ and $A$ and $B$ are clopen sets in $(G, \rho)$ and  $(M, \tau)$, respectively}
\end{definition}
The set of all clopen protoconcepts is denoted as $\mathfrak{P}^{T}(\mathbb{K}^{T})$ and the set of all clopen semiconcepts is denoted as $\mathfrak{H}^{T}(\mathbb{K}^{T})$.

\begin{theorem}
  \label{chartrep}  {\rm Let $\mathbb{K}^{T}(\textbf{D}):=((\mathcal{F}_{pr}(\textbf{D}),\mathcal{T}),(\mathcal{I}_{pr}(\textbf{D}),\mathcal{J}),\Delta)$ be a CTS based on the dBa {\bf D} and let $A\subseteq \mathcal{F}_{pr}(\textbf{D})$ and $B\subseteq \mathcal{I}_{pr}(\textbf{D})$. Then, we have 
    \begin{enumerate}
        \item If $\textbf{D}$ is fully contextual then $(A, B)$ is a clopen protoconcept if and only if $A=F_{x}$ and $B=I_{x}$ for some $x\in D$

         \item If $\textbf{D}$ is pure then $(A, B)$ is a clopen semiconcept if and only if $A=F_{x}$ and $B=I_{x}$ for some $x\in D$
    \end{enumerate}}
\end{theorem}
\begin{proof}
(1)   Let $(A, B)$ be a clopen  protoconcept of $\mathbb{K}^{T}(\textbf{D})$. Now we consider the clopen semiconcept $(A, A^{\prime})$ and $(B^{\prime}, B )$ then $(A, A^{\prime})_{\sqcup}=(A, A^{\prime})\sqcup (A, A^{\prime})=(A^{\prime\prime}, A^{\prime})=(B^{\prime}, B^{\prime\prime})=(B^{\prime}, B)\sqcap (B^{\prime}, B)=(B^{\prime}, B)_{\sqcap}$. As $A$ is clopen in  $(\mathcal{F}_{pr}(\textbf{D}),\mathcal{T})$  and $B$ is clopen in $(\mathcal{I}_{pr}(\textbf{D}),\mathcal{J})$, By Lemma \ref{clopenset}, $A=F_{x}$ and $B=I_{y}$ for $x, y\in D$. By Lemma \ref{derivation}(1), $(A, A^{\prime})=(F_{x}, I_{x\sqcap x})=(F_{x\sqcap x}, I_{x\sqcap x})=(F_{a}, I_{a})$, where $a=x\sqcap x$. y Lemma \ref{derivation} (2)$(B^{\prime}, B)=(F_{y\sqcup y}, I_{y})=(F_{y\sqcup y}, I_{y\sqcup y})=(F_{b}, I_{b})$, where $b=y\sqcup y$.

Now $(F_{a},I_{a})_{\sqcup}=(F_{a\sqcup a}, I_{a\sqcup a})$ and similarly, we can show that $(F_{b}, I_{b})_{\sqcap}=(F_{b\sqcap b}, I_{b\sqcap b})$.  So $(F_{a\sqcup a}, I_{a\sqcup a})=(F_{b\sqcap b}, I_{b\sqcap b})$, which implies that $a\sqcup a=b\sqcap b$. As $\textbf{D}$ is fully contextual there exist a $c\in D$, $c\sqcap c=a$ and $c\sqcup c=b$ which implies that $A=F_{a}=F_{c\sqcap c}=F_{c}$ and $B=I_{b}=I_{c\sqcup c}=I_{c}$. Hence $(A, B)=(F_{c}, I_{c})$.

Conversely, $F_{x}^{\prime\prime}=I_{x\sqcap x}^{\prime}=F_{(x\sqcap x)\sqcup (x\sqcap x)}=F_{x\sqcup x}=I_{x}^{\prime}$ and hence $(F_{x}, I_{x})$ is a clopen protoconcept.

\noindent (2) As $\textbf{D}$ is pure then $x\sqcap x=x$  or $x\sqcup x=x$ for $x\in D$. If $x\sqcap x=x$ then $(F_{x}, I_{x})=(F_{x}, F^{\prime}_{x})$ and if $x\sqcup x=x$ then $(F_{x}, I_{x})=(I^{\prime}_{x}, I_{x})$.

If the semiconcept is $(A, A^{\prime})$ then we can write it as $(A, A^{\prime})=(F_{x}, F_{x}^{\prime})=(F_{x}, I_{x\sqcap x})=(F_{x\sqcap x}, I_{x\sqcap x})$. If $(B^{\prime}, B)$ then we can write $(B^{\prime}, B)=(F_{x\sqcup x}, I_{x\sqcup x})$.
\end{proof}

\begin{corollary}
\label{new:repfulpuredba}    {\rm
    
    \noindent \begin{enumerate}
        \item For a fully contextual dBa $\textbf{D}$, the set $\mathfrak{P}^{T}(\mathbb{K}^{T}(\textbf{D}))$ forms a fully contextual dBa and is isomorphic to $\textbf{D}$.
        \item For a  pure dBa $\textbf{D}$, the set $\mathfrak{H}^{T}(\mathbb{K}^{T}(\textbf{D}))$ forms a pure dBa and is isomorphic to $\textbf{D}$.
    \end{enumerate}
    
     }
\end{corollary}
\begin{proof}
   (1) By Theorem \ref{chartrep}, the set $\mathfrak{P}^{T}(\mathbb{K}^{T}(\textbf{D}))$ coincides with the set  $\mathcal{D}$. By Corollary \ref{constalgeb}, $\mathfrak{P}^{T}(\mathbb{K}^{T}(\textbf{D}))$ forms a dBa. As fully contextual dBa is contextual, by Theorem \ref{new:rep} $\textbf{D}$ is isomorphic to $\mathfrak{P}^{T}(\mathbb{K}^{T}(\textbf{D}))$. Now it is remain to show that  the algebra is fully contextual. Let $(F_{a}, I_{a})$ and $(F_{b}, I_{b})$ belongs to $\mathfrak{P}^{T}(\mathbb{K}^{T}(\textbf{D}))$ such that $(F_{a}, I_{a})_{\sqcap}=(F_{a}, I_{a})$ and $(F_{b}, I_{b})_{\sqcup}=(F_{b}, I_{b})$ and $(F_{a}, I_{a})_{\sqcup}=(F_{b}, I_{b})_{\sqcap}$. So $(F_{a\sqcap a},I_{a\sqcap a})=(F_{a}, I_{a})$, $(F_{b\sqcup b},I_{b\sqcup b})=(F_{b}, I_{b})$ and $(F_{a\sqcup a}, I_{a\sqcup a})=(F_{b\sqcap b}, I_{b\sqcap b})$, which implies that $a\sqcap a=a$, $b\sqcup b=b$ and $a\sqcup a=b\sqcap b$. As $\textbf{D}$ is fully contextual there exist a $c\in D$, $c\sqcap c=a$ and $c\sqcup c=b$ which implies that $(F_{c}, I_{c})_{\sqcap}=(F_{c\sqcap c}, I_{c\sqcap c})=(F_{a}, I_{a})$ and $(F_{c}, I_{c})_{\sqcup}=(F_{c\sqcup c}, I_{c\sqcup c})=(F_{b}, I_{b})$.\\
   (2) The proof is similar to the proof of (1) and hence we omit the proof.
\end{proof}

Corollary~\ref{new:repfulpuredba} is an alternative to Theorems~\ref{iso-fullycxt dBa} and~\ref{RTDBA}. 
It provides a representation theorem for fully and pure dBas in terms of classical protoconcepts and semiconcepts of a context. 
The proof presented here is simpler than the one in~\cite{MR4566932}. 
Moreover,  the representation results presented in Theorems~\ref{iso-fullycxt dBa} and~\ref{RTDBA} can be obtained from Corollary~\ref{new:repfulpuredba}, and vice versa.

Now recall Theorem~\ref{connectingthm}(1), which states that the necessity, sufficiency, and derivation operators are dual in nature, and that each can be translated into the others. Keeping this observation in mind, we present a dual formulation of Definition~\ref{def:CTSCR topcon} in what follows.

\begin{definition}
   \label{tr-ctscr} {\rm A CTS \(\mathbb{K}^{T} := ((G, \rho), (M, \tau), R)\)  is called {\it translated-CTSCR}  if \(\mathbb{K}^{T} := ((G, \rho), (M, \tau), R^{c})\) is a CTSCR.}
    \end{definition}
Recall the operations $\sqcap$, $\sqcup$, $\neg$, $\lrcorner$, $\top$, and $\bot$ defined on the set $\mathfrak{P}(\mathbb{K})$ in Section~\ref{intro}. We then have the following result.
\begin{theorem}
   \label{protofully} {\rm Let \(\mathbb{K}^{T} := ((G, \rho), (M, \tau), R)\) be a translated-CTSCR. Then we have the following.
    \noindent \begin{enumerate}
 			\item $\underline{\mathfrak{P}}^{T}(\mathbb{K}^{T}):=(\mathfrak{P}^{T}(\mathbb{K}^{T}),\sqcup,\sqcap,\neg,\lrcorner,\top,\bot)$ is  a fully contextual dBa. Moreover, $\underline{\mathfrak{P}}^{T}(\mathbb{K}^{T})$  is isomorphic to $\underline{\mathfrak{R}}^{T}(\mathbb{K}^{T})$.
 			\item  $\underline{\mathfrak{H}}^{T}(\mathbb{K}^{T}):=(\mathfrak{H}^{T}(\mathbb{K}^{T}),\sqcup,\sqcap,\neg,\lrcorner,\top,\bot)$ is a subalgebra of  $\underline{\mathfrak{P}}^{T}(\mathbb{K}^{T})$ and is  a pure dBa. Moreover $\underline{\mathfrak{H}}^{T}(\mathbb{K}^{T})$ is isomorphic to $\underline{\mathfrak{S}}^{T}(\mathbb{K}^{T})$.
 	\end{enumerate}}
\end{theorem}
\begin{proof}
    (1) It is a routine check that the set $\mathfrak{P}^{T}(\mathbb{K}^{T})$ is closed under the operations 
$\sqcap$, $\sqcup$, $\neg$, and $\lrcorner$. 
For example, let $(A,B),(C,D)\in \mathfrak{P}^{T}(\mathbb{K}^{T})$. Then
$(A,B)\sqcap (C,D) = (A\cap C,(A\cap C)^{\prime})$.
Since $A$ and $C$ are clopen, $A\cap C$ is clopen. 
By Theorem~\ref{connectingthm}(1),
$(A\cap C)_{R}^{\prime} = (A\cap C)_{R^{c}}^{\boxg}$.
By Definition~\ref{tr-ctscr}, \(\mathbb{K}^{T} := ((G, \rho), (M, \tau), R)\) is a CTSCR, which implies that 
$(A\cap C)_{R^{c}}^{\boxg}$ is clopen. Hence, $(A\cap C)_{R}^{\prime}$ is clopen which implies that $(A,B)\sqcap (C,D)\in \mathfrak{P}^{T}(\mathbb{K}^{T})$.  The algebra $\underline{\mathfrak{P}}^{T}(\mathbb{K}^{T})$ is a contextual dBa as it is a subalgera of  $\underline{\mathfrak{P}}(\mathbb{K})$. To show it is a fully contextual dBas, let $(A, B), (C, D)\in \mathfrak{P}^{T}(\mathbb{K}^{T})$ such that $(A, B)_{\sqcap}=(A, A^{\prime})=(A, B)$ and $(C, D)_{\sqcup}=(D^{\prime}, D)=(C, D)$ and $(A, B)_{\sqcap}=(A^{\prime\prime}, A^{\prime})=(C, D)_{\sqcup}=(D^{\prime}, D^{\prime\prime})$ then $(A, D)\in \mathfrak{P}^{T}(\mathbb{K}^{T})$. Moreover $(A, D)_{\sqcap}=(A, B)$ and $(A, D)_{\sqcup}=(C, D)$. Hence it is a fully contextual dBa.

\end{proof}

\begin{note}
  Representation theorems  \ref{iso-fullycxt dBa} and  \ref{RTDBA} can be obtained as corollaries of Theorem \ref{protofully}
   and  Corollary \ref{new:repfulpuredba}. 
\end{note}



\section{On the Logics for Contextual and Pure dBas}\label{see:logic}
\subsection{Redefine the Logic for Contextual dBas}
 The proofs of soundness and completeness follow exactly the same pattern as those for \textbf{CDBL} in \cite{howlader2021dbalogic}. Therefore, we only present the proof system here and indicate that several axioms of \textbf{CDBL} are derivable within our new system. For details on \textbf{CDBL}, we refer the reader to \cite{howlader2021dbalogic}.

The alphabet of the logic $\textbf{L}$ consists of a countable set $P=\{p,q, \ldots\}$ of propositional variables propositional constants  $\bot,\top$, and logical connectives $\sqcup,\sqcap,\neg,\lrcorner$. Then, the set $Fm$ of formulas  is  defined inductively for $p\in P$ by
\[\varphi::= p\;\mid\;\neg\varphi\;\mid\;\lrcorner\varphi\;\mid\varphi\sqcap\varphi\;\mid\;\varphi\sqcup\varphi\;\mid\;\top\;\mid\bot.\]

\noindent The {\it axioms}  of $\textbf{L}$ are given by the following schema. 
\begin{center}
\[\varphi \Rightarrow \varphi;\quad
  \varphi \sqcap \psi \Rightarrow \varphi;\quad
  \varphi \sqcap \psi \Rightarrow \psi\]
  \end{center}
  \begin{center}
\[\varphi \Rightarrow \varphi \sqcup \psi;\quad
  \psi \Rightarrow \varphi \sqcup \psi;\quad
  \neg(\varphi \sqcap \varphi) \Rightarrow \neg\varphi;\quad
  \lrcorner\varphi \Rightarrow \lrcorner(\varphi \sqcup \varphi)\]
  \end{center}
  \begin{center}
  \[\varphi \sqcap \neg\varphi \Leftrightarrow \bot;\quad
  \top \Leftrightarrow \varphi \sqcup \lrcorner\varphi;\quad
  \neg\neg(\varphi \sqcap \psi) \Leftrightarrow \varphi \sqcap \psi;\quad
  \lrcorner\lrcorner(\varphi \sqcup \psi) \Leftrightarrow \varphi \sqcup \psi\]
  \end{center}
  \begin{center}
  \[\varphi \sqcap \varphi \Rightarrow \varphi \sqcap (\varphi \sqcup \psi);\quad
  \varphi \sqcup (\varphi \sqcap \psi) \Rightarrow \varphi \sqcup \varphi;\quad
  \varphi \sqcap (\psi \vee \theta) \Leftrightarrow (\varphi \sqcap \psi) \vee (\varphi \sqcap \theta);\]
  \end{center}
  \begin{center}
  \[\varphi \sqcup (\psi \wedge \theta) \Leftrightarrow (\varphi \sqcup \psi) \wedge (\varphi \sqcup \theta);\quad
  (\varphi \sqcup \varphi) \sqcap (\varphi \sqcup \varphi) \Leftrightarrow (\varphi \sqcap \varphi) \sqcup (\varphi \sqcap \varphi)\]
  \end{center}

\noindent {\it Rules of inference of {\bf L}} are as follows. 
    \begin{center}
        \[\infer[(Cut)]{  \varphi\Rightarrow\theta }{ \varphi\Rightarrow\psi &   \psi\Rightarrow\theta}\]
    \end{center}

\begin{center}
    \[\infer[(\sqcap_{R})]{\varphi\sqcap \theta\Rightarrow \psi\sqcap \theta}{\varphi\Rightarrow \psi} \quad
	\infer[(\sqcap_{L})]{ \theta\sqcap\varphi \Rightarrow \theta\sqcap \psi}{ \varphi\Rightarrow \psi}\quad
	\infer[(\sqcup_{R})]{  \varphi\sqcup \theta\Rightarrow \psi\sqcup \theta}{  \varphi\Rightarrow \psi}\]
\end{center}

\begin{center}
    \[\infer[(\sqcup_{L})]{  \theta\sqcup\varphi \Rightarrow\theta \sqcup \psi}{  \varphi\Rightarrow \psi}\quad\infer[(\neg)]{  \neg\psi\Rightarrow \neg\varphi}{  \varphi\Rightarrow \psi}\quad\infer[(\lrcorner)]{ \lrcorner\psi\Rightarrow \lrcorner\varphi}{ \varphi\Rightarrow \psi}\]
\end{center}
\vspace{.5cm}
\begin{center}
	    {\small $\infer[\hspace*{-4pt}(\sqsubseteq)]{ \varphi\Rightarrow\psi }{  \varphi\sqcap\psi\Rightarrow\varphi\sqcap\varphi &   \varphi\sqcap\varphi\Rightarrow\varphi\sqcap\psi &  \varphi\sqcup\psi\Rightarrow\psi\sqcup\psi &   \psi\sqcup\psi\Rightarrow\varphi\sqcup\psi}$}
\end{center}
Derivability is defined in the standard manner: an sequent $\varphi\Rightarrow\psi$ is {\it derivable} (or
{\it provable}) in \textbf{L}, if there exists a finite sequence of sequents $S_{1},. . ., S_{m}$
such that $S_{m}$ is $\varphi\Rightarrow\psi$ and for all $k \in\{1,. . .,m\}$ either $S_{k}$ is an axiom or $S_{k}$ is obtained by applying  rules of  \textbf{L} to elements from $\{S_{1} ,. . . , S_{k-1}\}$.
\begin{lemma}\label{lma:cpdl}
   \label{axiomcdbl1} {\rm For $\varphi,\psi,\theta\in Fm$, the following are provable in  \textbf{L}. \\
    $\begin{array}{ll}
    (1a)~\varphi\sqcap\psi\Rightarrow (\varphi\sqcap\psi)\sqcap(\varphi\sqcap\psi)& (1b)~(\varphi\sqcup\psi)\sqcup(\varphi\sqcup\psi)\Rightarrow\varphi\sqcup \psi.
    \end{array}$
    }
\end{lemma}

\begin{proof}
We give the proof of (1a); the proof of (1b) is dual.
\[\begin{array}{ll}
\text{1}\quad 
   \varphi \sqcap \psi \;\Rightarrow\; \neg\neg(\varphi \sqcap \psi)& \text{Axiom}\\  
\text{2}\quad
   \neg\!\left((\varphi \sqcap \psi)\sqcap(\varphi \sqcap \psi)\right)
   \;\Rightarrow\;
   \neg(\varphi \sqcap \psi)& \text{Axiom}\\
\text{3}\quad
   \neg\neg(\varphi \sqcap \psi)
   \;\Rightarrow\;
   \neg\neg\!\left((\varphi \sqcap \psi)\sqcap(\varphi \sqcap \psi)\right) & 2,(\neg)\\ 
\text{4}\quad
   \neg\neg\!\left((\varphi \sqcap \psi)\sqcap(\varphi \sqcap \psi)\right)
   \;\Rightarrow\;
   (\varphi \sqcap \psi)\sqcap(\varphi \sqcap \psi)& \text{Axiom}\\
\text{5}\quad \neg\neg(\varphi \sqcap \psi)
   \;\Rightarrow\;
   (\varphi \sqcap \psi)\sqcap(\varphi \sqcap \psi)& 3,4, \text{(Cut)}\\
\text{6}\quad \varphi \sqcap \psi
   \;\Rightarrow\;
   (\varphi \sqcap \psi)\sqcap(\varphi \sqcap \psi)& 1,5, \text{(Cut)}
\end{array}\]
\end{proof}
\begin{theorem} 
	\label{thempdbl}
	{\rm 
		\label{other axiom of DBA}
		For $\varphi,\psi\in Fm$, the following are provable in  \textbf{L}. 
		
		$	\begin{array}{ll}

		(1a)~ (\varphi\sqcap \psi)\Leftrightarrow(\psi\sqcap \varphi).&
		(1b)~ \varphi\sqcup \psi\dashv\Rightarrow \psi\sqcup \varphi.\\
		(2a)~ \neg \Leftrightarrow\neg (\varphi\sqcap \varphi).&
		(2b)~ \lrcorner(\varphi\sqcup \varphi)\Rightarrow\lrcorner \varphi.\\
		
		(3a)~ \varphi\sqcap(\varphi\sqcup \psi)\Rightarrow (\varphi\sqcap \varphi).&
		(3b)~ \varphi\sqcup \varphi\Rightarrow \varphi\sqcup (\varphi\sqcap \psi).\\
		\end{array}$
	}
\end{theorem}
\begin{proof} For $i\in \{1, 2,3\}$, we will give the proof of ($i$a); that of ($i$b) is dual.\\
 (1a) 
$\begin{aligned}[t]
\text{1} &\quad \varphi\sqcap\psi \;\Rightarrow\; \varphi& \mbox{Axiom}\\
\text{2} &\quad \varphi\sqcap\psi \;\Rightarrow\; \psi& \mbox{Axiom}\\
\text{3} &\quad (\varphi \sqcap \psi) \sqcap (\varphi \sqcap \psi) \;\Rightarrow\; (\psi\sqcap\varphi) & \sqcap_{L}, \sqcap_{R}, \mbox{and (Cut)}\\
\text{4} &\quad (\varphi \sqcap \psi) \;\Rightarrow\; (\varphi \sqcap \psi) \sqcap (\varphi \sqcap \psi) & \mbox{Lemma~\ref{lma:cpdl}(1a)}\\
\text{5} &\quad (\varphi \sqcap \psi) \;\Rightarrow\; (\psi\sqcap\varphi)& \mbox{4, 3, and (Cut)}
 \end{aligned}$l

\noindent	Interchanging $\varphi$ and $\psi$ in the above, we get  $(\psi\sqcap\varphi)\Rightarrow(\varphi\sqcap\psi)$. 

\noindent (2a)
$\begin{aligned}[t]
 \text{1} &\quad \varphi \sqcap \varphi \;\Rightarrow\; \varphi &\mbox{Axiom}\\
\text{2} &\quad \neg \varphi \;\Rightarrow\; \neg (\varphi \sqcap \varphi)& (\neg).
\end{aligned}$
		
\noindent (3a) 
$\begin{aligned}[t]
   \text{1} &\quad \varphi \sqcap (\varphi \sqcup \psi) \;\Rightarrow\; \varphi& \mbox{Axiom}\\
\text{2} &\quad \varphi \sqcap (\varphi \sqcup \psi) \;\Rightarrow\; \varphi& \mbox{Axiom}\\
\text{3} &\quad  
   (\varphi \sqcap (\varphi \sqcup \psi)) \sqcap (\varphi \sqcap (\varphi \sqcup \psi)) \;\Rightarrow\; \varphi \sqcap \varphi&  \sqcap_{L}, \sqcap_{R}, \text{and (Cut)}\\
\text{4} &\quad 
   \varphi \sqcap (\varphi \sqcup \psi) \;\Rightarrow\; (\varphi \sqcap (\varphi \sqcup \psi)) \sqcap (\varphi \sqcap (\varphi \sqcup \psi))& \mbox{Lemma!\ref{lma:cpdl}(1a)}\\
\text{5} &\quad  
   \varphi \sqcap (\varphi \sqcup \psi) \;\Rightarrow\; \varphi \sqcap \varphi& \mbox{4, 3, (Cut)} 
\end{aligned}
$

\end{proof}
\begin{definition}
    {\rm Let $\textbf{D}$ be a contextual D-core  algebra. A sequent $\varphi\Rightarrow\psi$ is said to be satisfied by a homomorphism $h:Fm\rightarrow\textbf{D}$ if $h(\varphi)\sqsubseteq h(\psi)$. In addition, $\varphi\Rightarrow\psi$ is called true in $\textbf{D}$ if it is satisfied by every homomorphism $h:Fm\rightarrow\textbf{D}$.   $\varphi\Rightarrow\psi$ is called valid if it is true in any contextual D-core  algebra.}
\end{definition}
\begin{theorem}
   \label{sundcomp} {\rm A sequent $\varphi\Rightarrow\psi$ is provable if and only it is valid.}
\end{theorem}

Next theorem ensure that despite these simplifications, the two proof systems $\textbf{CDBL}$ and $\textbf{L}$ are equivalent, so the essential logical properties are preserved.
\begin{theorem}
    \label{axiomcdbl2}{\rm For $\varphi,\psi\in Fm$, the following are provable in  \textbf{L}. \\
    $\begin{array}{ll}
    (1a)~\varphi\sqcap \varphi\Rightarrow \varphi\sqcap(\varphi\vee \psi)&
	(1b)~\varphi\sqcup(\varphi\wedge \psi)\Rightarrow \varphi\sqcup \varphi\\
    (2a)~  \bot\Rightarrow \varphi &
	(2b)~ \varphi\Rightarrow \top\\
	(3a)~\neg\top\Rightarrow\bot &
	(3b)~\top\Rightarrow\lrcorner\bot\\
	(4a)~\neg \bot\Leftrightarrow \top\sqcap\top &
	(4b)~\lrcorner\top\Leftrightarrow\bot\sqcup\bot\\
    \end{array}$
    }
\end{theorem}
\begin{proof}
    {\rm  Items 1a and 1b follow from Theorem~\ref{thm:axiom6ad} together with Theorem~\ref{sundcomp}, while items 3a, 3b, 4a, and 4b are also consequences of Theorem~\ref{thm:firstindeax} and Theorem~\ref{sundcomp}. Items 2a and 2b, on the other hand, follow from Proposition~\ref{proaxiom6}(4a and 4b) and Proposition~\ref{prop:axion6}(1a, 3a, 1b, and 3b)  in combination with Theorem~\ref{sundcomp}.
 }
\end{proof}
\begin{note}
    {\rm Note that the proof system \(\textbf{L}\) is a simplified version of the proof system \(\mathbf{CDBL}\). In particular, \(\textbf{L}\) can be obtained from \(\mathbf{CDBL}\) by removing the redundant axioms listed in Lemma~\ref{axiomcdbl1} and Theorem~\ref{axiomcdbl2}. Moreover, Lemma~\ref{axiomcdbl1} and Theorem~\ref{axiomcdbl2} ensure that the two systems are equivalent.
}
\end{note}
\subsection{Redefine the Logic for pure dBas}\label{see:plogic}
In \cite{HOWLADER2023115}, a hypersequent calculus \textbf{PDBL} is developed for pure double Boolean algebras. The alphabet of \textbf{PDBL} consists of a countably infinite set \(\mathbf{OV} = \{p, q, r, \ldots\}\) of object variables, a countably infinite set \(\mathbf{PV} = \{P, Q, R, \ldots\}\) of property variables, the propositional constants \(\bot\) and \(\top\), and the logical connectives \(\sqcup, \sqcap, \neg,\) and \(\lrcorner\). The set \(Fm\) of formulas is defined inductively by 
\[\varphi::= P\;\mid\; p\;\mid\;\neg\varphi\;\mid\;\lrcorner\varphi\;\mid\varphi\sqcap\varphi\;\mid\;\varphi\sqcup\varphi\;\mid\;\top\;\mid\bot,\] where \(p \in \mathbf{OV}\) and \(P \in \mathbf{PV}\). An s-hypersequent in \textbf{PDBL} is a finite sequence of sequents of the form \(\varphi_{1} \Rightarrow \psi_{1} \mid \varphi_{2} \Rightarrow \psi_{2} \mid \ldots \mid \varphi_{n} \Rightarrow \psi_{n}\). Note that a sequent  $\varphi\Rightarrow \psi$  is a special case of a hypersequent when $n=1$.

The axioms closely resemble those of \textbf{CDBL} \cite{howlader2021dbalogic}, which is expected since every pure dBa is also a contextual one. In the case of \textbf{PDBL}, two additional axioms are introduced to handle the variables. 

\[p\sqcap p\Leftrightarrow p; \quad
 P\sqcup P\Leftrightarrow P\]
The syntactic differences between \textbf{PDBL} and \textbf{CDBL} are reflected primarily in the rules of inference.
\begin{center}
    \[\infer[(\sqcap_{R})]{B\mid\varphi\sqcap \theta\Rightarrow\psi\sqcap \theta\mid  C}{B\mid  \varphi\Rightarrow\psi\mid  C} \quad
	\infer[(\sqcap_{L})]{B\mid  \theta\sqcap\varphi \Rightarrow\theta\sqcap \psi\mid  C}{B\mid  \varphi\Rightarrow\psi\mid  C}\quad\infer[(\neg)]{B\mid  \neg\psi\Rightarrow\neg\varphi\mid  C}{B\mid  \varphi\Rightarrow\psi\mid  C}
	\]
\end{center}

\begin{center}
    \[\infer[(\sqcup_{L})]{B\mid  \theta\sqcup\varphi \Rightarrow\theta \sqcup \psi\mid  C}{B\mid  \varphi\Rightarrow\psi\mid  C}\quad\infer[(\sqcup_{R})]{B\mid  \varphi\sqcup \theta\Rightarrow\psi\sqcup \theta\mid  C}{B\mid  \varphi\Rightarrow\psi\mid  C}\quad	\infer[(\lrcorner)]{B\mid  \lrcorner\psi\Rightarrow\lrcorner\varphi\mid  C}{B\mid  \varphi\Rightarrow\psi\mid  C}\]
\end{center}

\begin{center}
    \[\infer[(Cut)]{B\mid  D\mid  \varphi\Rightarrow\theta\mid  C\mid  E}{B\mid  \varphi\Rightarrow\psi\mid  C & D\mid  \psi\Rightarrow\theta\mid  E}\quad\infer[(Sp)]{\varphi\Rightarrow\varphi\sqcap\varphi\mid  \varphi\sqcup\varphi\Rightarrow\varphi}{}\]
\end{center}

{\tiny\begin{center}
    \[\infer[(\sqsubseteq)]{B\mid  D\mid  F\mid  H\mid  \varphi\Rightarrow\psi\mid  C\mid  E\mid  G\mid  X}{B\mid  \varphi\sqcap\psi\Rightarrow\varphi\sqcap\varphi\mid  C & D\mid  \varphi\sqcap\varphi\Rightarrow\varphi\sqcap\psi\mid  E & F\mid  \varphi\sqcup\psi\Rightarrow\psi\sqcup\psi\mid  G & H\mid  \psi\sqcup\psi\Rightarrow\varphi\sqcup\psi\mid  X}\]
\end{center}}
{\it External rules of inference}:
\begin{center}
    \[\infer[\mbox{(EC)}]{B\mid D\mid  C}{B\mid  D\mid  D\mid  C}\quad \infer[\mbox{(EE)}]{B\mid  E\mid  D\mid  C}{B\mid  D\mid  E\mid  C}\quad\infer[\mbox{(EW)}]{B\mid C}{B}\]
\end{center}
Derivability is defined in the standard way. The calculus is shown to be sound and complete with respect to the class of pure dBas. For the detailed proofs, we refer the reader to \cite{HOWLADER2023115}.

As in the previous section, we can obtain a simplified hypersequent calculus \textbf{HL} for pure dBas that is equivalent to \textbf{PDBL}. This is done by taking the axioms of \textbf{L} instead of all the axioms of \textbf{CDBL}. Since sequents are themselves hypersequents, Lemma~\ref{axiomcdbl1} and Theorem~\ref{axiomcdbl2} also hold in the new proof system \textbf{HL}. Consequently, the system \textbf{HL} is equivalent to \textbf{PDBL}. As \textbf{PDBL} is sound and complete with respect to pure dBas, the same holds for \textbf{HL}. The proofs of soundness and completeness are exactly the same as those for \textbf{PDBL}, and therefore we omit the details here.

\section{Conclusion}\label{see:concl}
 We provide the first reduced axiom system for double Boolean algebras, refine existing representation theorems, generalize the classical glued-sum construction, establish a Stone-type  topological representation theorem, and introduce simplified logical systems with complete proof theories. These results significantly advance the algebraic, topological, and logical foundations of double Boolean algebras and open new directions for research in Contextual Logic, Boolean Concept Logic, and related areas of algebraic logic.

\subsection*{Acknowledgment}
This research was funded by the SNF Grant IZSEZO\textunderscore233403/1:
\emph{Axiomatizations of Double Boolean Algebras}.



\section*{Declarations}

\subsection*{Ethical approval}
Not applicable. 

\subsection*{Competing interests} 
Not applicable. 
\subsection*{Authors' contributions} 
All authors contributed equally.

\subsection*{Availability of data and materials}
Not applicable. 
\smallskip





\bibliographystyle{plain}
\bibliography{name}

@article {MR4566932,
    AUTHOR = {Howlader, Prosenjit and Mohua Banerjee},
     TITLE = {Topological representation of double {B}oolean algebras},
   JOURNAL = {Algebra Universalis},
  FJOURNAL = {Algebra Universalis},
    VOLUME = {84},
      YEAR = {2023},
      PAGES = {Paper No. 15, 32},
      }

@article{kembang2023simple,
  title={Simple and sub-directly irreducible double Boolean algebras},
  author={Kembang, GT and Kwuida, L and Temgoua, ERA and Tenkeu, YLJ},
  journal={Preprint arXiv:2312.13686},
  year={2023}
}

@article{bmlpl,
    author = {P. Howlader and C.J. Liau},
    title = {On the Logical and Algebraic Aspects of Reasoning with Formal Contexts},
    year = {2025},
    address = {New York, NY, USA},
    volume = {26},
    number = {3},
    journal = {ACM Trans. Comput. Logic},
    articleno = {17}
    }

@misc{howladerdiscussion,
  title={A Discussion on Double Boolean Algebras Extended Abstract},
  author={Howlader, Prosenjit and Liau, Churn-Jung},
  journal={TACL 2024},
  year={2024},
  howpublished={\url{https://iiia.csic.es/tacl2024/abstracts/conference/contributed/TACL_2024_paper_115.pdf}}
}

@incollection{Lakhal2005,
author="L. Lakhal and G. Stumme",
editor="B. Ganter and G. Stumme  and R. Wille",
title="Efficient mining of association rules based on formal concept analysis",
bookTitle="Formal Concept Analysis: Foundations and Applications",
year="2005",
publisher="Springer",
address="Berlin, Heidelberg",
pages="180--195",
}

@book{birkhoff1940lattice,
  title={Lattice theory},
  author={Birkhoff, Garrett},
  volume={25},
  year={1940},
  publisher={American Mathematical Soc.}
}

@article{HOWLADER2023115,
title = {A non-distributive logic for semiconcepts and its modal extension with semantics based on {K}ripke contexts},
journal = {International Journal of Approximate Reasoning},
volume = {153},
pages = {115-143},
year = {2023},
author = {Prosenjit Howlader and Mohua Banerjee},

}

@article{howlader2021dbalogic,
title = {Kripke Contexts, Double {B}oolean Algebras with Operators and Corresponding Modal Systems},
journal = {Journal of Logic, Language and Information},
volume = {32},
pages = {117-146},
year = {2023},
author = {Prosenjit Howlader and Mohua Banerjee}}

@article{pawlak1982rough,
  title="Rough sets",
  author="Pawlak, Zdzis{\l}aw",
  journal="International Journal of Computer and Information Sciences",
  volume="11",
  number="5",
  pages="341--356",
  year="1982",
 }

@InProceedings{howlader3,
author = "Howlader, Prosenjit and Banerjee, Mohua",
title = "Remarks on Prime Ideal and Representation Theorems for Double {B}oolean Algebras",
editor="Francisco J. Valverde-Albacete and others",
booktitle="CLA 2020",
pages = "83-94",
publisher="CEUR Workshop Proceedings",
year = "2020",
}

@InProceedings{howlader2020,
author="Howlader, Prosenjit
and Banerjee, Mohua",
title="Object Oriented Protoconcepts and Logics for Double and Pure Double {B}oolean Algebras",
editor="Bello, Rafael  and others",
booktitle="Int. Joint Conf. on Rough Sets",
pages="308--323",
publisher="Springer, Cham",
year="2020",
}

@inproceedings{howlader2018algebras,
  author="Howlader, Prosenjit and
   Banerjee, Mohua",
  title="Algebras from Semiconcepts in Rough Set Theory",
  editor=" Hung Son Nguyen and  others",
  booktitle="International Joint Conference on Rough Sets",
  pages="440--454",
  publisher="Springer, Cham",
  year="2018",
}

@incollection{wille1982restructuring,
 author={Wille, Rudolf},
  title={Restructuring lattice theory: an approach based on hierarchies of concepts},
  editor={Rival, I.},
  booktitle={ Ordered Sets. NATO Advanced Study Institutes Series (Series C — Mathematical and Physical Sciences)},
  pages={445--470},
  publisher={Springer, Dordrecht},
  year={1982},
}

@book{guide2006infinite,
 author={Aliprantis, Charalambos D. and  Border, Kim C.},
title={Infinite Dimensional Analysis: A Hitchhiker's Guide},
publisher={Springer, Berlin },
  year={2006}
  }

@inproceedings{duntsch2002modal,
 author={D{\"u}ntsch, Ivo   and Gediga, G{\"u}nther},
  title={Modal-style operators in qualitative data analysis},
  booktitle={Proceedings of the 2002 {IEEE} Int. Conf. on Data Mining},
  editor={Vipin, Kumar and others},
  pages={155--162},
  publisher={{IEEE} Computer Society},
  year={2002},
}

@InProceedings{wille,
author="Wille, Rudolf",
title="Boolean Concept Logic",
editor="Ganter, Bernhard and others",
booktitle="Conceptual Structures: Logical, Linguistic, and Computational Issues",
pages="317--331",
publisher="Springer, Berlin",
year="2000",
}

@article{kwuida2007prime,
  author={Kwuida, L{\'e}onard},
  title={Prime ideal theorem for double {B}oolean algebras},
  journal={Discussiones Mathematicae-General Algebra and Applications},
  volume={27},
  number={2},
  pages={263--275},
  year={2007},
  
}

@article{BALBIANI2012260,
author = "Philippe Balbiani",
title = "Deciding the word problem in pure double {B}oolean algebras",
journal = "Journal of Applied Logic",
volume = "10",
number = "3",
pages = "260 - 273",
year = "2012",

}

\end{document}